\documentclass[12pt,a4paper,leqno]{amsart}
\usepackage{amssymb}
\usepackage{amsmath}
\usepackage{amsthm}
\usepackage{indentfirst}

\newtheorem{lemma}{Lemma}
\newtheorem{theorema}{Theorem}

\newtheorem{theoremb}{Theorem}

\newtheorem*{claim*}{Claim}

\theoremstyle{remark}
\newtheorem*{remark*}{Remark}

\numberwithin{equation}{section}
\numberwithin{lemma}{section}
\newenvironment{lista}
{\begin{list}{} {\setlength{\itemsep}{0 ex plus 0.1 ex}
\setlength{\labelwidth}{1 cm} \setlength{\topsep}{1 ex}
\setlength{\partopsep}{0 pt}}} {\end{list}}
\begin{document}
\date{}
\title[Derivatives of Meromorphic Functions]{Derivatives of Meromorphic Functions with Multiple Zeros and Elliptic Functions}

 \author[Pai Yang]{Pai Yang}
 \email{yangpai@cuit.edu.cn}


\author[Shahar Nevo]{Shahar Nevo
 }

 \address{Department of Mathematics,
 Bar-Ilan University, 52900 Ramat-Gan, Israel}
  \email{nevosh@macs.biu.ac.il}

  \author[Xuecheng Pang]{Xuecheng Pang
  }
  \address{Department of Mathematics,
East China Normal University,  Shanghai 200062, P.~R. China}
\email{xcpang@euler.math.ecnu.edu.cn}

\begin{abstract}
 Let $f$ be a nonconstant meromorphic function in the plane and $h$ be a nonconstant elliptic function.
 We show that if all zeros of $f$ are multiple except finitely many and $T(r,h)=o\{T(r,f)\}$ as $r\rightarrow\infty$,
 then $f'=h$ has infinitely many solutions (including poles).
\end{abstract}



\maketitle
\section{Introduction}
Recall that an elliptic function \cite{Elliptic199079} is a
meromorphic function $h$ defined on $\mathbb{C}$ for which there
exist two non-zero complex numbers $\omega_1$ and $\omega_2$ with
$\omega_1/\omega_2$ not real such that $h(z+\omega_1) =
h(z+\omega_2)=h(z)$   for all $z$ on $\mathbb{C}$.

We use the following notation. For $z_0\in \mathbb{C}$ and $r>0$,
$\Delta(z_0,r)=\{z:|z-z_0|<r\}$,
$\Delta'(z_0,r)=\{z:0<|z-z_0|<r\}$, $\Delta=\Delta(0,1)$,
$\overline \Delta(0,r)=\{z:|z|\leq r\}$, and
$\Gamma(0,r)=\{z:|z|=r\}$. For $f$ meromorphic in a domain  $D$,
we denote
$$S(D,f)=\frac{1}{\pi}\iint_D[f^\#(z)]^2dxdy, \qquad  S(r,f)=S(\Delta(0,r),f).$$
Here
$$f^\#(z)=\frac{|f'(z)|}{1+|f(z)|^2}$$
denotes the spherical derivative.
We use the Ahlfors-Shimizu form of the Nevanlinna characteristic function, given by
$$T_0(r,f)=\int_0^r\frac{S(t,f)}{t}dt.$$
Let $n(r,f)$ denote the number of poles of $f(z$) in $\Delta(0,r)$
(counting multiplicity). We write
$f_n\overset{\chi}{\Longrightarrow}f$ in $D$ to indicate that the
sequence $\{f_n\}$ converges to $f$ in the spherical metric
uniformly on compact subsets of $D$ and $f_n\Rightarrow f$ in $D$
if the convergence is in the Euclidean metric.

In 1959, Hayman \cite{hayman195970} proved the following seminal result,
which has come to be known as Hayman's Alternative.
\begin{theorema}\label{thayman}
  Let $f$ be a transcendental meromorphic function on the complex plane $\mathbb{C}$. Then either
\begin{lista}
 \item [$(i)$] $f$ assumes each value $a\in \mathbb{C}$ infinitely often, or
 \item [$(ii)$] $f^{(k)}$ assumes each value $b\in \mathbb{C} \backslash\{0\}$ infinitely often for $k=1,2,\ldots$.
\end{lista}
\end{theorema}
Obviously, if $f$ is a transcendental meromorphic function such that $f\neq 0$ for all $z\in \mathbb{C}$,
then $f^{(k)}$ assumes each value $b\in \mathbb{C} \backslash\{0\}$ infinitely often for $k=1,2,\ldots$.
In recent years, it has become clear that, in many instances, the
condition $f\neq 0$ can be replaced by the assumption that all zeros
of $f$ have sufficiently high multiplicity. This announcement
concerns such extension of Theorem \ref{thayman}. We restrict our
attention to the case $k=1$.

Before stating our result, let us present several theorems that
show results already obtained in this subject.

\begin{theorema}\textup{\cite[Theorem 3]{wang199814}} \label{twanga}
Let $f$ be a transcendental meromorphic function on $\mathbb{C}$, all of whose zeros have multiplicity at least 3. Then $f'$ assumes each nonzero complex value infinitely often.
\end{theorema}
The analogue of Theorem \ref{twanga} with 3 replaced by 2 does
hold for functions of finite order.
\begin{theorema}\textup{\cite[Lemma 6]{wang199814}} \label{twangb}
Let $f$ be a transcendental meromorphic function of finite order on $\mathbb{C}$, all of whose zeros are multiple. Then $f'$ assumes every finite non-zero value infinitely often.
\end{theorema}
This is an instant corollary of Theorem \ref{thayman} and the following important result.
\begin{theorema}\textup{\cite[Theorem 3]{bergweiler199511}} \label{tbergweilera}
Let $f$ be a transcendental meromorphic function of finite order on $\mathbb{C}$ with an infinite number of multiple zeros. Then $f'$ assumes every finite non-zero value infinitely often.
\end{theorema}
Indeed, if in Theorem \ref{twangb}, $f$ vanishes only finitely often,
then $f'$ must take on every nonzero value infinitely often by Theorem \ref{thayman};
Otherwise, Theorem \ref{tbergweilera} implies the same conclusion.

Bergweiler and Eremenko gave a counter example in \cite{bergweiler199511}
which shows Theorem \ref{tbergweilera} is not true in general for functions of infinite order.

In 2006, Shahar Nevo, Xuecheng Pang and Lawrence Zalcman gave the
following extension of Theorem \ref{twanga} and Theorem
\ref{twangb}.
\begin{theorema}\textup{\cite[Theorem 1]{pang20082}}\label{trational}
Let $f$ be a transcendental meromorphic function on $\mathbb{C}$, all but
finitely many of whose zeros are multiple, and let $R\not\equiv0 $
be a rational function. Then $f'-R$ has infinitely many zeros.
\end{theorema}
In this paper, we continue to study the value distribution of the
derivative of meromorphic functions with multiple zeros. We have the
following theorem.
\begin{theoremb}\label{thm}
 Let $f$ be a nonconstant meromorphic function on $\mathbb C$ and $h$ be a nonconstant elliptic function.
 Then if all zeros of $f$ are multiple except finitely many and $T_0(r,h)=o\{T_0(r,f)\}$ as $r\rightarrow\infty$,
 then $f'=h$ has infinitely many solutions (including the possibility of infinitely many common poles of $f$ and $h$).
\end{theoremb}

\begin{remark*}
$T(r,f)$ denotes the usual Nevanlinna characteristic function.
Since $T(r,f)-T_0(r,f)$ is bounded as a function of $r,$ one can
replace $T_0(r,f)$ with $T(r,f)$ in Theorem~\ref{thm}.
\end{remark*}

\section{Preliminary results}
\begin{lemma}\textup{\cite[Theorem 1]{tsuji19502}}\textup{\cite[Theorem VI.21]{tsuji19752}}\label{lalforscovera}
Let $f $ be a meromorphic function in $\Delta,$ and let
$a_1,a_s,a_3$ be $3$ distinct complex numbres. Assume that the
number of zeros of $\prod\limits_{i=1}^3(f(z)-a_i)$ in $\Delta$ is
$\leq n$, where multiple zeros are counted only once, then
\begin{equation}\nonumber
  S(r,f)\leq n+\frac{A}{1-r},\quad 0\leq r <1,
\end{equation}
and $A>0$ is a constant, which depends on $a_1,a_2,a_3$ only.
\end{lemma}
\begin{lemma}\label{lelliptic}
If $h(z)$ is a nonconstant elliptic function with primitive
periods $\omega_1$, $\omega_2$, where $\omega_1/\omega_2$ not
real, then $$T_0(r,h) = Ar^2(1 + o(1))\quad \mbox{as}\quad
r\rightarrow \infty,$$ where $A>0$ is a constant.
\end{lemma}
This follows from \cite[Corollary 2]{Bank1984981}.
\begin{lemma}\textup{\cite[Theorem 3]{schwick199732}}\textup{\cite[Theorem 1]{Shahara20014}}\label{lschwick}
Let $\psi\not\equiv 0$ be a meromorphic function in a domain $D
\subset \mathbb{C}$ and $k\in \mathbb{N}$. Let $\mathcal F$ be a
family of meromorphic functions in $D$, such that $f$ and
$f^{(k)}-\psi$ have no zeros and $f$ and $\psi$ have no common
poles for each $f\in \mathcal F$. Then $\mathcal F$ is normal on
$D$.
\end{lemma}
\begin{lemma}\label{lrationlwang}
  Let $f$ be a nonconstant meromorphic function of finite order on $\mathbb{C}$, all of whose zeros are multiple.
  If $f'(z)\neq 1$ on $\mathbb{C}$, then
  $$f(z)=\frac{(z-a)^2}{z-b}$$
  for some $a$ and $b\,(\neq a)$ on $\mathbb{C}$.
\end{lemma}
  This follows from Lemma $6$ (with $j=1$ and $k=2$) and Lemma $8$ (with $k=1$) of \cite{wang199814}.
\begin{lemma}\textup{\cite[Lemma 2]{pang200032}} \label{zp}
Let $\mathcal F$ be a family of functions meromorphic in a domain
$D$, all of whose zeros have multiplicity at least $k$, and
suppose that there exists $A\geq1$ such that $|f^{(k)}(z)|\leq A$
whenever $f(z)=0$.
Then if $\mathcal F$ is not normal at $z_0$, there exist, for each $0\leq\alpha\leq k$,\\
$(a)$ points $z_n$, $z_n\rightarrow z_0$;\\
$(b)$ functions $f_n\in\mathcal F$; and\\
$(c)$ positive numbers $\rho_n\rightarrow 0$\\
such that
$\rho_n^{-\alpha}f_n(z_n+\rho_n\zeta)=g_n(\zeta)\overset{\chi}{\Longrightarrow}g(\zeta)$ on $\mathbb{C}$,
where $g$ is a nonconstant meromorphic function on $\mathbb{C}$, all of whose zeros have
multiplicity at least $k$, such that $g^\#(\zeta)\leq g^\#(0)=kA+1$.
\end{lemma}
\begin{lemma}\label{lnormala}\textup{\cite[Lemma 3.1]{pang200554}}
  Let $\{f_n\}$ be a family of functions meromorphic in $\Delta(z_0,r)$, all of whose zeros and poles are multiple; and let $\{b_n\}$ be a sequence of holomorphic functions in $\Delta(z_0,r)$ such that $b_n\Rightarrow b$ in $\Delta(z_0,r)$, where $b(z)\neq 0$,
  $z\in \Delta(z_0,r)$. If $f'_n(z)\neq b_n(z)$ for $z\in \Delta(z_0,r)$, then $\{f_n\}$ is normal in $\Delta(z_0,r)$.
\end{lemma}
\begin{lemma}\textup{\cite[Lemma 3.1]{pang2005158}}\label{guyongxing}
Let $\{f_n\}$ be a sequence of meromorphic functions in
$\Delta(z_0,r)$ and $\{\psi_n\}$ a sequence of holomorphic
functions in $\Delta(z_0,r)$ such that $\psi_n\Rightarrow\psi$,
where $\psi(z)\not\neq 0$ in $\Delta(z_0,r)$. If for each $n$,
$f_n(z)\neq 0$ and $f'_n(z)\neq\psi_n(z)$ for $z\in\Delta(z_0,r)$,
then $\{f_n\}$ is normal in $\Delta(z_0,r)$.
\end{lemma}
\begin{lemma}\textup{\cite[Lemma 3.2]{pang2005158}}\label{lquasiconezeroa}
Let $E$ be a (countable) discrete set in $\Delta(z_0,r)$ which has
no accumulation points in $\Delta(z_0,r)$ and let $\{\psi_n\}$ be
a sequence of holomorphic functions in $\Delta(z_0,r)$ such that
$\psi_n\Rightarrow\psi$ in $\Delta(z_0,r)$, where $\psi\neq
0,\infty$ in $\Delta(z_0,r)$. Let $\{f_n\}$ be a sequence of
functions meromorphic in $\Delta(z_0,r)$, all of whose zeros are
multiple, such that $f'_n(z)\neq\psi_n(z)$ for all $n$ and all
$z\in \Delta(z_0,r)$. Let $a_1\in E$ and suppose that
\begin{lista}
  \item [$(a)$] for some $a_1\in E$, no subsequence of $\{f_n\}$ is   normal at $a_1$;
  \item [$(b)$] there exists $\delta>0$ such that each $f_n$ has a single (multiple) zero in $\Delta(a_1,\delta)$; and
  \item [$(c)$] $f_n(z)\overset{\chi}{\Longrightarrow}f(z)$ in $\Delta(z_0,r)\backslash E$.
\end{lista}
Then
\begin{lista}
  \item [$(d)$] there exists $\eta_0>0$ such that for each $0<\eta<\eta_0$, $f_n$ has a single simple pole in $\Delta(a_1,\eta)$ for sufficiently large $n$; and
  \item [$(e)$] $f(z)=\int^z_{a_1}\psi(\zeta)d\zeta$.
\end{lista}
\end{lemma}
\begin{lemma}\textup{\cite[Lemma 3.7]{pang2005158}}\label{lquasictwozeroa}
Let $\{f_n\}$ be a sequence of functions meromorphic on
$\Delta(z_0,r)$, all of whose zeros are multiple, and let $\psi$
be a non-vanishing holomorphic function in $\Delta(z_0,r)$.
Suppose that
\begin{lista}
  \item [$(a)$] $\{f_n\}$ is quasinormal in $\Delta(z_0,r)$;
  \item [$(b)$] $f'_n(z)\neq\psi(z)$ for $z\in \Delta(z_0,r)$ and $n=1,2,3,\cdots$;
  \item [$(c)$] no subsequence of $\{f_n(z)\}$ is normal at $0$.
\end{lista}
Then there exists $0<\delta<r$ such that $f_n$ has only a single
(multiple) zero in $\Delta(0,\delta)$ for sufficiently large $n$.
\end{lemma}
\begin{lemma}\label{lquasinormalityshahar}
Let $\{f_n\}$ be a family of meromorphic functions on the plane
domain $D$, all whose zeros are multiple, If there exists a
holomorphic function $\varphi$ univalent in $D$ such that
$f'_n(z)\neq \varphi'(z)$ for each $n\in\mathbb{N}$ and $z\in D$,
Then $\{f_n\}$ is quasinormal of order $1$ in $D$.
\end{lemma}
This follows from Theorem $1$ (with $k=1$ and $\mathcal F=\{f_n\}$) of \cite{Shahar2007101}.
\begin{lemma}\textup{\cite[Lemma 12]{changjianming2008338}}\label{lrationala}
Let $R$ be a nonconstant rational function satisfying $R'\neq 0$ on $\mathbb{C}$.
Then $R(z)=az+b$ or $R(z)=\frac{a}{(z+c)^n}+b$, where $n\in \mathbb{N}$ and $a\neq 0,b,c\in\mathbb{C}$.
\end{lemma}
\section{Auxiliary lemmas}
\begin{lemma}\label{lrationalb}
Let $k\geq 2$ be an integer and let $R$ be a rational function on
$\mathbb{C}$.
Suppose that $R'(z)\neq\frac{1}{z^{k}}$. Then $R$ is a constant
function.
\end{lemma}
\begin{proof}
Since $R'(z)\neq\frac{1}{z^{k}}$, we have $R(0)\neq \infty$ and hence
$$ R'(z)-\frac{1}{z^{k}}=\big(R(z)-\frac{1}{1-k}\frac{1}{z^{k-1}}\big)'\neq 0.$$
 By Lemma \ref{lrationala},
 either $R(z)=az+b+\frac{1}{1-k}\frac{1}{z^{k-1}}$, or
 $R(z)=\frac{a}{(z+c)^n}+b+\frac{1}{1-k}\frac{1}{z^{k-1}}$.
 If  $R(z)=az+b+\frac{1}{1-k}\frac{1}{z^{k-1}}$, then we have
 $R(0)=\infty$. A contradiction.
 Thus, $R(z)=\frac{a}{(z+c)^n}+b+\frac{1}{1-k}\frac{1}{z^{k-1}}$. Since
 $R(0)\neq \infty$, we have that $c=0$,  $n=k-1$, $a=\frac{1}{k-1}$ and  $R(z)\equiv b$.
\end{proof}
\begin{lemma}\label{lrationalc}
Let $k$ be a positive integer and  let $R$ be a rational function
satisfying $R'(z)\neq z^k$ on $\mathbb{C}$. If all zeros of $R$
are  multiple, then
\begin{equation}\label{equationrational}
  R(z)=\frac{\prod\limits_{i=1}^{n+k+1}(z-\alpha_i)}{(k+1)(z-\beta)^n},
\end{equation}
 where $n$ is a nonnegative integer $\beta\in \mathbb{C}$ and  $\alpha_i(\neq0,\beta),$ for $1\le i\le n+k+1$ and every zero $\alpha_i$ is counted
 due to multiplicity.
\end{lemma}
\begin{proof}
 Obviously, $(R(z)-\frac{z^{k+1}}{k+1})'\neq 0$ and since all the zeros of $R$ are multiple, then $R(z)-\frac{z^{k+1}}{k+1}$
 is a nonconstant rational function. By Lemma~\ref{lrationala},  $R(z)=\frac{z^{k+1}}{k+1}+az+b$ or $R(z)=\frac{z^{k+1}}{k+1}+\frac{a}{(z+c)^n}+b$,
 where $n\in \mathbb{N}$ and $a\neq 0,b,c\in\mathbb{C}$.
 Let us consider separately these two cases.

{\bf Case 1:} $R(z)=\frac{z^{k+1}}{k+1}+az+b$.

Let $z_0$ be a zero of $R(z)$, then $R(z_0)=0$ and hence
$R'(z_0)=0$, i.e.,
\begin{align}\label{equationrationala}
  &R(z_0)=\frac{1}{k+1}z_0^{k+1}+az_0+b=0,\\
  &R'(z_0)=z_0^k+a=0.\nonumber
\end{align}
It follows that $z_0=-\frac{b(k+1)}{ak}$ and hence
\begin{equation}\label{equationrationalb}
  R(z)=\frac{1}{k+1}\left(z+\frac{b(k+1)}{ak}\right)^{k+1}.
\end{equation}
Comparing the coefficients of (\ref{equationrationala}) and (\ref{equationrationalb}),
we have $k=1$ and hence $R(z)=\frac{(z-\alpha)^2}{2}$, where $\alpha=-\frac{2b}{a}$.
Now $R(z)$ have the form of (\ref{equationrational}) with $\alpha_1=\alpha_2=\alpha$, $k=1$ and $n=0$.

 {\bf Case 2:} $R(z)=\frac{1}{k+1}z^{k+1}+\frac{a}{(z+c)^n}+b$.

 Let $\{\alpha_1, \alpha_2, \cdots, \alpha_{n+k+1}\}$ be the set of the zeros of  $R(z)$ and $\beta=-c$.
 Obviously, $R(z)$ have the form of (\ref{equationrational}).

 Since $R'(z)\neq z^k$,   \ $R'(0)\ne0$, and since all zeros of $R$ are  multiple,     $R(0)\neq 0$.
   Hence $\alpha_i\neq 0$ for $1\le i\le n+k+1$.
\end{proof}
\begin{lemma}\label{lnonormalconvert}
Let $\{f_n(z)\}$ be a family of meromorphic functions in $\Delta$,
all of whose zeros are multiple. Let $F_n(z)=z^kf_n(z)$, where
$k\neq 0$ is an integer. Let $\{b_n(z)\}$ be a sequence of
holomorphic functions in $\Delta$ such that $b_n(z)\Rightarrow
b(z)$ in $\Delta$, where $ b(z)\neq 0$ is a holomorphic function
in $\Delta$. Suppose that for each $n$, $f'_n(z)\neq z^{-k}b_n(z)$
and $\{f_n(z)\}$ is not normal at $0$ and normal in $\Delta'$.
Then $\{F_n(z)\}$ is   normal in $\Delta'$ but not normal at $0$.
\end{lemma}
\begin{proof}
It is obvious that $\{F_n(z)\}$ is normal in $\Delta'$. Suppose
that $\{F_n(z)\}$ is normal at $0$.

{\bf Case 1:} \textit{$k$ is a positive integer.}

Since $f'_n(z)\neq z^{-k}b_n(z)$ and $b_n(z)\Rightarrow b(z)$ on
$\Delta$, and $b(0)\ne0$ we have, for sufficiently large $n$,
$f'_n(0)\neq \infty$ and $f_n(0)\neq \infty$. Obviously, we have
$F_n(0)=0$. Without loss of generality, we can suppose that for
all $n$, $f_n(0)\neq \infty$ and $F_n(0)=0$. Since $\{F_n(z)\}$ is
normal at $0$ and $F_n(0)=0$, there exists $\delta>0$ such that
for all $n$, $|F_n(z)|\leq 1$ in $\Delta(0,\delta)$. Thus
$f_n(0)\neq \infty$ in $\Delta'(0,\delta)$, and hence $f_n$ is
holomorphic in $\Delta(0,\delta)$. Now, we have
$$|f_n(z)|=\frac{|F_n(z)|}{|z|^k}\leq \left(\frac{1}{\delta}\right)^k,\quad z\in \Gamma(0,\delta).$$
By the maximum principle, this holds throughout $\Delta(0,\delta)$.
It follows that $\{f_n\}$ is normal. A contradiction.

{\bf Case 2:}  \textit{$k$ is a negative integer.}

Let $m=-k$. Since $f'_n(z)\neq z^mb_n(z)$ and $b_n(z)\Rightarrow
b(z)$ in $\Delta$, we have, for sufficiently large $n$,
$f'_n(0)\neq 0$ and hence $f_n(0)\neq 0$. Obviously, we have
$F_n(0)=\infty$. Without loss of generality, we can suppose that
for all $n$, $f_n(0)\neq 0$ and $F_n(0)=\infty$. Since
$\{F_n(z)\}$ is normal at $0$ and $F_n(0)=\infty$, there exists
$\delta>0$ such that for all $n$, $|F_n(z)|\geq 1$ on
$\Delta(0,\delta)$. Thus $f_n(z)\neq 0$ in $\Delta'(0,\delta)$,
and hence $\frac{1}{f_n}$ is holomorphic in $\Delta(0,\delta)$.
 Now, we have
 $$\left|\frac{1}{f_n(z)} \right|=\left|\frac{1}{F_n(z)}\frac{1}{z^m} \right|\leq \frac{1}{\delta^m},\quad z\in \Gamma(0,\delta).$$
 By the maximum principle, this holds throughout $\Delta(0,\delta)$.
 It follows that $\{f_n\}$ is normal. A contradiction.
\end{proof}
\begin{lemma}\label{lhurwitz}
Let $\{f_n\}$ and $\{\psi_n\}$ be two sequences of functions
meromorphic on the plane domain $D$. Let $f(z)$ and $\psi(z)$ be
two meromorphic functions in $D$. Suppose that
\begin{lista}
 \item [$(a)$] $f_n(z) \overset{\chi}{\Longrightarrow} f(z)$ in $D$ and $\psi_n(z) \overset{\chi}{\Longrightarrow} \psi(z)$ in $D$;
 \item [$(b)$] $f'_n(z)\neq \psi_n(z)$ in $D$.
\end{lista}
Then, either $f'(z)\equiv \psi(z)$ in $D$, or $f'(z)\neq \psi(z)$
in $D$.
\end{lemma}
\begin{proof}
Suppose that there exists $z_0\in D$ such that
$f'(z_0)=\psi(z_0)$.

Let us separate into two cases.

\textbf{Case 1:} $\psi(z_0)\ne\infty.$

There exists $r>0$ such that $f,\psi,f_n,\psi_n$ are analytic in
$\Delta(z_0,r)$ for large enough $n.$ We have that
$f_n'-\psi_n\Rightarrow f'-\psi$ in $\Delta(z_0,r)$, and by
condition (b) and Hurwitz's Theorem $f'-\psi\equiv 0$ in
$\Delta(z_0,r)$ and so in all of $D.$ In this case, the first
possibility occurs.

\textbf{Case 2:} $\psi(z_0)=\infty.$

Suppose by negation that $f'-\psi$ is not constantly zero. Thus,
for small enough $r$, $\frac{1}{f'-\psi}$ is analytic in
$\Delta'(z_0,r)$, and since $\frac{1}{f_n'-\psi_n}$ are analytic
in $D$ and $\frac{1}{f_n'-\psi_n}\Rightarrow \frac{1}{f'-\psi}$ in
$\Delta'(z_0,r)$, we get by the maximum principle that
$\left\{\frac{1}{f_n'-\psi_n}\right\}_{n=1}^\infty$ is a uniformly
convergent Cauchy sequence of analytic functions in
$\Delta(z_0,r)$ which thus converges (uniformly) to an analytic
function there. Thus $\frac{1}{f'-\psi}$ is extended analytically
to $\Delta(z_0,r)$. Let $m_1,m_2\ge 1$ be the orders of the pole
of $f$ and $\psi$ at $z_0,$ respectively. A pole of $f'-\psi$ at
$z_0,$ if it exists, is of order at most $\max\{m_1+1,m_2\}.$ This
is also the maximal multiplicity of zero of $\frac{1}{f'-\psi}$ at
$z_0.$ (Observe that it is possible that $f'(z_0)-\psi(z_0)$ is a
finite value, but this value cannot be $0$ since
$\frac{1}{f'-\psi}$ is analytic at $z_0).$ Now, for every
$0<\delta<r$,  $f_n$ has, for large enough $n,$ $m_1$ poles
(counting multiplicities) in $\Delta(z_0,\delta),$ and $\psi_n$
has $m_2$ poles in $\Delta(z_0,\delta)$ (counting multiplicities).
Since $f_n$ and $\psi_n$ have no common poles, then
$\frac{1}{f_n'-\psi_n}$ has at least $m_1+m_2+1$ zeros in
$\Delta(z_0,\delta).$ But $m_1+m_2+1\ge \max\{m+1,m_2\}$ and this
is  a contradiction. Thus $f'-\psi\equiv0.$ So also in the case of
a common pole, we have that $f'\equiv\psi$ in $D$ and the lemma is
proved.
\end{proof}
\begin{lemma}\label{ldnormal}
  Let $\{f_n\}$ be a family of meromorphic functions in $\Delta$, all of whose zeros are multiple.
  Let $\{b_n\}$ be a sequence of holomorphic functions in $\Delta$ such that $b_n\Rightarrow 1$ in $\Delta$,
  and let $k$ be a positive integer.   Suppose that
 \begin{lista}
 \item [$(a)$] $f'_n(z)\neq z^kb_n(z)$ in $\Delta$;
 \item [$(b)$] there exists points $z_n\rightarrow 0$ such that $f_n(z_n)=0$;
 \item [$(c)$] $f_n(z)\overset{\chi}{\Longrightarrow}f(z)$ in $\Delta'$, where $f(z)$ is a meromorphic function in $\Delta'$.
\end{lista}
  Then $f(z)=\frac{z^{k+1}+c}{k+1}$ in $\Delta'$, where $c$ is a constant.
\end{lemma}
\begin{proof}
  Let $F_n(z)=\frac{f'_n(z)}{z^k}$.
  Since $b_n\Rightarrow 1$ in $\Delta$ and $(a)$,
  we have for sufficiently large $n$, $f'_n(0)\neq 0$ and hence $F_n(0)=\infty$.
  Without loss of generality, we may assume that for all $n$, $f'_n(0)\neq 0$ and $F_n(0)=\infty$.
  Since all zeros of $\{f_n(z)\}$ are multiple, $f_n(0)\neq 0$.
  Hence we have $z_n\neq 0$ and $F_n(z_n)=0$.

  We claim that  $\{F_n(z)\}$ is not normal at $0$ and hence $\{\frac{f'_n(z)}{z^kb_n(z)}-1\}$ is also not normal at $0$.
 Indeed, since $F_n(z_n)=0$ and $F_n(0)=\infty$, the family $\{F_n(\zeta)\}$ is not equicontinuous at $0$ and hence cannot be normal at $0$.

  By $(a)$ and $(c)$, we have
  $$0\neq\frac{f'_n(z)}{z^kb_n(z)}-1\Rightarrow \frac{f'(z)}{z^k}-1,\quad z\in \Delta'-\{f^{-1}(\infty)\}.$$
By Hurwitz's Theorem, either $\frac{f'(z)}{z^k}-1\equiv 0$ in
$\Delta'$, or $\frac{f'(z)}{z^k}-1\neq 0$ in $\Delta'$. Suppose
that $\frac{f'(z)}{z^k}-1\neq 0$ in $\Delta'$. Since $f(z)$ is a
meromorphic function, then there exists $\delta>0$ such that
$f(z)$ has no poles on $\Gamma(0,\delta)$ and $f'_n(z)$ converges
uniformly to $f'(\zeta)$ on $\Gamma(0,\delta)$. Now, we have
$$\infty\neq\frac{1}{\frac{f'_n(z)}{z^kb_n(z)}-1}\Rightarrow \frac{1}{\frac{f'(z)}{z^k}-1},\quad z\in\Gamma(0,\delta).$$
As in the proof of Lemma \ref{lhurwitz}, since the function in the
left hand side is holomorphic, we have by the maximum principle
that this holds throughout $\Delta(0,\delta)$. So
$\{\frac{f'_n(z)}{z^kb_n(z)}-1\}$ is normal at $0$. A
contradiction. Thus, $\frac{f'(z)}{z^k}-1\equiv 0$ in $\Delta'$.
Then $f(z)=\frac{z^{k+1}+c}{k+1}$ in $\Delta'$, where $c$ is a
constant.
\end{proof}
\begin{lemma}\label{lalforscoverb}
Let $\{f_n\}$ be a sequence of functions meromorphic in
$\Delta(z_0,r)$. Suppose there exists $M_1>0$ such that for each
$n$, $n(r,\frac{1}{f_n})<M_1$. Suppose that
$f_n\overset{\chi}{\Longrightarrow}f$ in $\Delta'(z_0,r)$, where
$f$ is a nonconstant meromorphic function or $f\equiv\infty$ in
$\Delta'(z_0,r)$. Then there exists $M_2>0$ such that, for
sufficiently large $n$,
$$S(\frac{r}{2}, f_n)<M_2.$$
\end{lemma}
\begin{proof}
 Without loss of generality, we assume that $r=1$ and so $\Delta(z_0,r)=\Delta.$ Obviously, $\frac{1} {f_n}-1\overset{\chi}{\Longrightarrow}\frac{1}{f}-1$ in $\Delta'$
  and $\frac{1}{f}-1\not\equiv 0,\infty$ in $\Delta'$.
  Then there exists $\frac{1}{2}<r<1$ such that $\frac{1}{f}-1$ has no poles and zeros in $\Gamma(0,r)$.
  Obviously, for sufficiently large $n$,
  \begin{align}\nonumber
    n\left(r,\frac{1}{f_n-1}\right)-n\left(r,\frac{1}{f_n}\right)&=n\left(r,\frac{1}{\frac{1}{f_n}-1}\right)-n\left(r,\frac{1}{f_n}-1\right)\\
    &=\frac{1}{2\pi i}\int_{\Gamma(0,\,r)}\frac{(\frac{1}{f_n}-1)'}{\frac{1}{f_n}-1}dz\rightarrow   \frac{1}{2\pi i}\int_{\Gamma(0,\,r)}\frac{(\frac{1}{f}-1)'}{\frac{1}{f}-1}dz. \nonumber
  \end{align}
  Set $M_3=\frac{1}{2\pi i}\int_{\Gamma(0,\,r)}\frac{(\frac{1}{f}-1)'}{\frac{1}{f}-1}dz+n(r,\frac{1}{f_n})+1$. For sufficiently large $n$,
 \begin{align*}\left(\frac{1}{2},\frac{1}{f_n-1}\right)&\leq
 n\left(r,\frac{1}{f_n-1}\right)
   =n\left(r,\frac{1}{f_n}\right)+\int_{\Gamma(0,r)}\frac{(\frac{1}{f}-1)'}{\frac{1}{f}-1}dz+\varepsilon_n\\ &<M_1+\int_{\Gamma(0,r)}
  \frac{(\frac{1}{f}-1)}{\frac{1}{f}-1}dz=M_3-1<M_3.\end{align*}
  (Here $\varepsilon_n
\to0,$ but since the other terms are integers, then
$\varepsilon_n=0$ for large enough $n.)$

  Obviously, $\frac{1} {f_n}-\frac{1}{2}\overset{\chi}{\Longrightarrow}\frac{1}{f}-\frac{1}{2}$ in $\Delta'$
  and $\frac{1}{f}-\frac{1}{2}\not\equiv 0,\infty$ in $\Delta'$.
  Then there exists $\frac{1}{2}<t<1$ such that $\frac{1}{f}-\frac{1}{2}$ has no poles and zeros in $\Gamma(0,t)$.
  Clearly, for sufficiently large $n$,
  \begin{align}\nonumber
    n(t,\frac{1}{f_n-2})-n\left(t,\frac{1}{f_n}\right)&=n\left(t,\frac{1}{\frac{1}{f_n}-\frac{1}{2}}\right
    )-n\left(t,\frac{1}{f_n}-\frac{1}{2}\right)\\
    &=\frac{1}{2\pi i}\int_{\Gamma(0,\,t)}\frac{(\frac{1}{f_n}-\frac{1}{2})'}{\frac{1}{f_n}-\frac{1}{2}}dz\rightarrow   \frac{1}{2\pi i}\int_{\Gamma(0,\,t)}\frac{(\frac{1}{f}-\frac{1}{2})'}{\frac{1}{f}-\frac{1}{2}}dz.\nonumber
  \end{align}
  Set $M_4=\frac{1}{2\pi i}\int_{\Gamma(0,\,t)}\frac{(\frac{1}{f}-\frac{1}{2})'}{\frac{1}{f}-\frac{1}{2}}dz+n(t,\frac{1}{f_n})+1$.
  Similarly to the previous paragraph, we have, for sufficiently large~$n$,
$$n\left(\frac{1}{2},\frac{1}{f_n-2}\right)\leq n\left(t,\frac{1}{f_n-2}\right)<M_4-1.$$

 Let $M_2=M_1+M_3+M_4+2A$. By Lemma  \ref{lalforscovera}, for sufficiently large $n$,
$$S(\frac{1}{2}, f_n)\leq
n(\frac{1}{2},\frac{1}{f_n})+n(\frac{1}{2},\frac{1}{f_n-1})+n(\frac{1}{2},\frac{1}{f_n-2})+2A<M_2.$$
\end{proof}
\begin{lemma}\label{lunivalent}
Let $\Psi(z)$ be a holomorphic univalent function  on
$\Delta(0,R)$ and $\{\Psi_n(z)\}$ be a sequence of holomorphic
functions in $\Delta(0,R)$ such that $\Psi_n(z)\Rightarrow
\Psi(z)$ in $\Delta(0,R)$. Then for each $r\in (0,R)$, we have,
for sufficiently large $n$,
\begin{lista}
  \item [$(a)$] $\Psi_n(z)$ is a holomorphic univalent function  in $\Delta(0,r)$;
  \item [$(b)$] there exists  $\delta_1>0$ such that $\Delta(\Psi(0),\delta_1)\subset \Psi(\Delta(0,r))$ and  $\Delta(\Psi(0),\delta_1)\subset \Psi_n(\Delta(0,r))$;
   \item [$(c)$]  there exists $\delta_2>0$ such that $\Psi(\Delta(0,\delta_2))\subset \Delta(\Psi(0),\delta_1)$ and $\Psi_n(\Delta(0,\delta_2))\subset \Delta(\Psi(0),\delta_1)$.
\end{lista}
\end{lemma}
\begin{proof}
Let $r<r_1<R$. Suppose that there exists a subsequence of
$\{\Psi_n(z)\}$ (that we continue to call $\{\Psi_n(z)\}$) which
is not univalent in $\Delta(0,r)$. Then there exist distinct
complex numbers $z_{n,1}$ and $z_{n,2}$ in $\Delta(0,r)$ such that
$\Psi_n(z_{n,1})=\Psi_n(z_{n,2})$. Obviously,
$n(r_1,\frac{1}{\Psi_n(z)-\Psi_n(z_{n,1})})~\geq~2$. Without loss
of generality, we may suppose that $z_{n,1}\rightarrow z_1$ and
$z_{n,2}\rightarrow z_2$ as $n\rightarrow\infty$. Obviously
$|z_1|\leq r$, $|z_2|\leq r$. Since $\Psi_n(z)\Rightarrow \Psi(z)$
in $\Delta(0,R)$, we have
$$\Psi(z_1)=\lim\limits_{n\rightarrow\infty}\Psi_n(z_{n,1})=\lim\limits_{n\rightarrow\infty}\Psi_n(z_{n,2})=\Psi(z_2).$$
Hence we have $z_1=z_2$.

 Obviously, $\Psi_n(z)-\Psi_n(z_{n,1})\Rightarrow \Psi(z)-\Psi(z_1)$ in $\Delta(0,R)$
 and $\Psi(z)-\Psi(z_1)\neq 0$ on $\Gamma(0,r_1)$.
 By the argument principle, for sufficiently large~$n$,
  $$n\left(r_1,\frac{1}{\Psi_n(z)-\Psi_n(z_{n,1})}\right)=n\left(r_1,\frac{1}{\Psi(z)-\Psi(z_1)}\right)=1.$$
  A contradiction occurs and hence $(a)$ holds.

 Obviously, there exists $\delta_1$ such that  $\Delta(\Psi(0),2\delta_1) \subset \Psi(\Delta(0,r))$,
 and hence $\Delta(\Psi(0),\delta_1) \subset \Psi(\Delta(0,r))$.
 For each $w\in\Delta(\Psi(0),\delta_1)$, we have $|\Psi(z)-w|>\delta_1$ on $\Gamma(0,r)$,
 and we have for sufficiently large $n$,
 $$|(\Psi_n(z)-w)-(\Psi(z)-w)|=|\Psi_n(z)-\Psi(z)|<\delta_1,\quad z\in\Gamma(0,r).$$
 By the argument principle, for sufficiently large $n$,
  $$n\left(r,\frac{1}{\Psi_n(z)-w}\right)=n\left(r,\frac{1}{\Psi(z)-w}\right)=1.$$
 This shows that for each $w\in \Delta(\Psi(0),\delta_1)$, there exists $z_0\in \Delta(0,r)$
 such that $\Psi_n(z_0)=w$. Hence $(b)$ holds.

Obviously, there exists $\delta_2$ such that
$\Psi(\Delta(0,\delta_2))\subset\Delta(\Psi(0),\delta_1/2)$, and
hence $\Psi(\Delta(0,\delta_2))\subset\Delta(\Psi(0),\delta_1)$.
Since $\Psi_n(z)\Rightarrow\Psi(z)$ in $\Delta(0,R)$, we have that
for sufficiently large $n$, $|\Psi_n(z)-\Psi(z)|<\delta_1/2$ in
$\Delta(0,\delta_2)$ and hence
$$|\Psi_n(z)-\Psi(0)|\leq |\Psi_n(z)-\Psi(z)|+|\Psi(z)-\Psi(0)|<\delta_1/2+\delta_1/2=\delta_1$$
in $\Delta(0,\delta_2)$. Hence $(c)$ holds.
\end{proof}
\begin{lemma}\label{lquasichiefa}
Let $\{f_n\}$ be a sequence of meromorphic functions in
$\Delta(z_0,r)$, all of whose zeros are multiple, and let
$\{\psi_n\}$ be a sequence of meromorphic functions in
$\Delta(z_0,r)$ such that $\psi_n\Rightarrow\psi$ in
$\Delta(z_0,r)$, where $\psi$ is a non-vanishing holomorphic
function in $\Delta(z_0,r)$. Let $E$ be a (countable) discrete set
in $\Delta(z_0,r)$ which has no accumulation points in
$\Delta(z_0,r)$. Suppose that
\begin{lista}
  \item [$(a)$] $f_n(z)\overset{\chi}{\Longrightarrow}f(z)$ in $\Delta(z_0,r)\backslash E$;
  \item [$(b)$] for some $a_1\in E$, no subsequence of $\{f_n\}$ is normal at $a_1$;
  \item [$(c)$] for all $n\in \mathbb{N}$, $f'_n(z)\neq\psi_n(z)$ in $\Delta (z_0,r)$.
\end{lista}
Then
\begin{lista}
    \item [$(d)$]  There exists $r\geq 0$ such that for sufficiently large $n$, $f_n$ has a single zero $z_{n,1}$ of order $2$
     and a single pole $z_{n,2}$ of order $1$ in $\Delta(a_1,r)$, where $z_{n,i}\rightarrow a_1$ as $n\rightarrow\infty$, $i=1,2$;
  \item [$(e)$] $f(z)=\int^z_{a_1}\psi(\zeta)d\zeta$.
\end{lista}
\end{lemma}
\begin{proof}
Without loss of generality, we may suppose that $a_1=0$. Set
$\Psi(z)=\int_{\zeta=0}^z\psi(\zeta) d\zeta$. There exists
$\frac{1}{2}>\delta>0$ such that $\Psi(z)$ is a univalent function
in $\Delta(0,2\delta)$. Since $\psi_n(z)\Rightarrow\psi(z)$ on
$\Delta$, without loss of generality, we may suppose that for all
$n\in \mathbb{N}$, $\psi_n$ is a holomorphic function in
$\Delta(0,2\delta)$.

  Let $\Psi_n(z)=\int_{\zeta=0}^z\psi_n(\zeta) d \zeta$, $z\in \Delta(0,2\delta)$.
  Obviously, we have that $\Psi_n(z)\Rightarrow \Psi(z)$ in $\Delta(0,2\delta)$.
  By Lemma \ref{lunivalent}, we have for sufficiently large $n$,
  \begin{lista}
  \item [$(a)$] $\Psi_n(z)$ is a holomorphic univalent function  in $\Delta(0,\delta)$;
  \item [$(b)$] there exists $\delta_1>0$ such that $\Delta(\Psi(0),\delta_1)\subset \Psi(\Delta(0,\delta))$ and $\Delta(\Psi(0),\delta_1)\subset \Psi_n(\Delta(0,\delta))$;
   \item [$(c)$] there exists $\delta_2>0$ such that $\Psi(\Delta(0,\delta_2))\subset \Delta(\Psi(0),\delta_1)$ and $\Psi_n(\Delta(0,\delta_2))\subset \Delta(\Psi(0),\delta_1)$.
\end{lista}
 For convenience, we suppose that for all $n$, $(a)$, $(b)$ and $(c)$ hold.

 Now we consider $\Psi(z)$ and $\Psi_n(z)$ only in $\Delta(0,\delta)$.
 Let $F_n(w)=f_n(\Psi^{-1}_n(w))$, $w\in \Delta(\Psi(0),\delta_1)$.
 Then $f_n(z)=F_n(\Psi_n(z))$ in $\Delta(0,\delta_2)$.

 We claim that no subsequence of $\{F_n(w)\}$ is normal at $\Psi(0)$.
 Otherwise, suppose that $\{F_n(w)\}$ is normal at $\Psi(0)$ (without loss of generality, we call also  the subsequence
 $\{F_n(w)\})$. Since
 $\Psi_n(z)\Rightarrow \Psi(z)$ in $\Delta(0,2\delta)$, $\{F_n(\Psi_n(z))\}$ is normal at $0$.
 i.e., $\{f_n(z)\}$ is normal at $0$, a contradiction.

 Now we have
  $$F'_n(w)= f'_n(\Psi^{-1}_n(w))\frac{1}{\Psi'_n(z)}=f'_n(z)\frac{1}{\psi_n(z)}\neq \psi_n(z)\frac{1}{\psi_n(z)}=1 $$
 in $\Delta(\Psi(0),\delta_1)$. Obviously, all the zeros of $F_n(w)$ in $\Delta(\Psi(0),\delta_1)$ are multiple.
 Clearly, $\{F_n\}_{n=1}^\infty$ is normal in $\Delta'(\psi(0),\delta_1).$ Then by Lemma~\ref{lquasictwozeroa}, there exists
$0<\delta_3<\delta_1$ such that $F_n(w)$ has only a single
(multiple) zero in $\Delta(\Psi(0),\delta_3)$ for sufficiently
large $n$.
 Do as in Lemma \ref{lunivalent}, we have that for sufficiently large $n$, there exists $\delta_4>0$ such that $\Psi_n(\Delta(0,\delta_4))\subset \Delta(\Psi(0),\delta_3)$. Therefore, $f_n(z)=F_n(\Psi_n(z))$ has at most a single (multiple) zero  in $\Delta(0,\delta_4)$.

 We claim that for sufficiently large $n$, $f_n$ has a single zero $z_{n,1}$ of order 2 in $\Delta(0,\delta_4)$, where $z_{n,1}\rightarrow 0$ as $n\rightarrow\infty$.
It suffices to prove that each subsequence of $\{f_n\}$ has a
subsequence $\{f_m\}$ such that, for sufficiently large $m$, $f_m$
has a single zero $z_{m,1}$ of order 2 in $\Delta(0,\delta_4)$,
where $z_{m,1}\rightarrow 0$ as $m\rightarrow\infty$. Suppose that
we have a subsequence of $\{f_n\}$, which (to avoid complication
in notion) we again call $\{f_n\}$. Since $\{f_n\}$ is not normal
at $0$, it follows from Lemma \ref{zp} that we can extract a
subsequence $\{f_m\}$ of $\{f_n\}$, points $z_m\rightarrow 0$, and
positive numbers $\rho_m\rightarrow 0$ such that
$$g_m(\zeta)=\frac{f_m(z_m+\rho_m\zeta)}{\rho_m}\overset{\chi}{\Longrightarrow}g(\zeta),\quad \zeta\in \mathbb{C},$$
where $g$ is a nonconstant meromorphic function of finite order on $\mathbb{C}$, all of whose zeros are multiple.
Now we have
$$0\neq f'_m(z_m+\rho_m\zeta)-\psi_m(z_m+\rho_m\zeta)\Rightarrow g'(\zeta)-\psi(0),\quad \zeta\in \mathbb{C}-\{g^{-1}(\infty)\}.$$
By Hurwitz's Theorem, either $g'(\zeta)\neq \psi(0)$ on
$\mathbb{C}$, or $g'(\zeta)\equiv \psi(0)$  on $\mathbb{C}$. In
the latter case, $g(\zeta)=\psi(0)\zeta+c$, which contradicts the
fact that all zeros of $g$ are multiple. Thus $g'(\zeta)\neq
\psi(0)$ on $\mathbb{C}$. By Lemma \ref{lrationlwang},
$g(\zeta)=\frac{\psi(0)(\zeta-a)^2}{\zeta-b}$ for distinct complex
numbers $a$ and $b$. It now follows from Hurwitz's Theorem that
there exist sequence $\zeta_{m,1}\rightarrow a$ such that, for
sufficiently large $m$, $g_m(\zeta_{m,1})=0$. Obviously,
$\zeta_{m,1}$ is a zero of order 2 of $g_m(\zeta)$. Set
$z_{m,1}=z_m+\rho_m\zeta_{m,1}$, we have $z_{m,1}\rightarrow 0$
and $f_m(z_{m,1})=0$. Obviously $z_{m,1}$ is a zero of order 2 of
$f_m(z)$ and since we have proved that there is at most one such
zero in $\Delta(0,\delta_4)$ (for sufficiently large $m$), it is
the only one, as required.

 By Lemma \ref{lquasiconezeroa},
 there exists $\eta>0$ such that $f_n$ has a single simple pole $z_{n,2}$ in $\Delta(0,\eta)$ for sufficiently large $n$.
 We claim that $z_{n,2}\rightarrow 0$ as $n\rightarrow 0$.
 Otherwise, there exist $0<\delta_5<\delta_4$ and a subsequence of $\{f_n\}$ (that we continue to call $\{f_n\}$)
 such that $f_n(z)\neq \infty$, $z\in \Delta(0,\delta_6)$.
 By Lemma \ref{lnormala}, $\{f_n\}$ is normal at $0$. A contradiction.
Set $r=\min\{\delta_4,\,\eta\}$. Hence $(d)$ holds.

 By Lemma \ref{lquasiconezeroa}, we have $f(z)=\int^z_{a_1}\psi(\zeta)d\zeta$. Hence $(e)$ holds.
\end{proof}
\begin{lemma}\label{lquasichiefb} (cf. \cite[Lemma
7]{Shahar2007101})
Let $\{f_n\}$ be a family of meromorphic functions on the plane
domain $D$, all of whose zeros are multiple, and let $\{\psi_n\}$
be a sequence of meromorphic functions in $D$ such that
$\psi_n\overset{\chi}{\Longrightarrow}\psi$ in $D$, where
$\psi(z)\not\equiv 0,\infty$ in $D$. If for each $n\in
\mathbb{N}$, $f'_n\neq \psi_n$ in $D$, then $\{f_n\}$ is
quasinormal in $D$.
\end{lemma}
\begin{proof}
 It suffices to show that $\{f_n\}$ is quasinormal in a neighborhood of each point of $D$.
 Let $p\in D$. There exists $t>0$ such that $\Delta(p,t)\subset D$
 and $\psi$ is holomorphic and does not vanish in $\Delta'(p,t)$.

 For each $q\in \Delta'(p,t)$,
 Let $\Psi(z)=\int_{\zeta=q}^z\psi(\zeta)d\zeta$ in $\Delta'(p,t)$.
 Since $\psi(z)$ is holomorphic and does not vanish in $\Delta'(p,t)$ and $\psi_n\overset{\chi}{\Longrightarrow}\psi$ in $D$,
 there exists $0<R<t-|p-q|$ such that for sufficiently large $n$,
 $\Psi(z)$ is a holomorphic univalent function  in $\Delta(q,R)$ and $\psi_n(z)$ is a holomorphic function in $\Delta(q,R)$.
 Let $\Psi_n(z)=\int_{\zeta=q}^z\psi_n(\zeta) d \zeta$ in $\Delta(q,R)$.
 Obviously, we have $\Psi_n(z)\Rightarrow\Psi(z)$ in $\Delta(q,R)$.

 Let $0<r<R$. By Lemma \ref{lunivalent}, we have, for sufficiently large $n$,
 \begin{lista}
  \item [$(a)$] $\Psi_n(z)$ is a holomorphic univalent function  in $\Delta(q,r)$;
  \item [$(b)$] there exists $\delta_1$ such that $\Delta(\Psi(q),\delta_1)\subset \Psi(\Delta(q,r))$ and $\Delta(\Psi(q),\delta_1)\subset \Psi_n(\Delta(q,r))$;
   \item [$(c)$] there exists $\delta_2$ such that $\Psi(\Delta(q,\delta_2))\subset \Delta(\Psi(q),\delta_1)$ and $\Psi_n(\Delta(q,\delta_2))\subset \Delta(\Psi(q),\delta_1)$.
\end{lista}

  For convenience in notation, we assume that for all $n$, $(a)$, $(b)$ and $(c)$ hold.

  Now we consider $\Psi(z)$ and $\Psi_n(z)$ only in $\Delta(q,r)$.
  Let $F_n(w)=f_n(\Psi^{-1}_n(w))$, $w\in \Delta(\Psi(q),\delta_1)$.
  Then $f_n(z)=F_n(\Psi_n(z))$ in $\Delta(q,\delta_2)$.

  Now we have
  $$F'_n(w)= f'_n(\Psi^{-1}_n(w))\frac{1}{\Psi'_n(z)}=f'_n(z)\frac{1}{\psi_n(z)}\neq \psi_n(z)\frac{1}{\psi_n(z)}=1 $$
  in $\Delta(\Psi(q),\delta_1)$. Obviously all the zeros of $F_n(w)$ are multiple in $\Delta(\Psi(q),\delta_1)$.
  By Lemma \ref{lquasinormalityshahar}, $\{F_n(w)\}$ is quasinormal in $\Delta(\Psi(q),\delta_1)$.
  Since $f_n(z)=F_n(\psi_n(z))$ in $\Delta(q,\delta_2)$ and $\psi_n(z)\Rightarrow \psi(z)$ in $\Delta(q,R)$,
 $\{f_n(z)\}$ is quasinormal in $\Delta(q,\delta_2)$ and hence quasinormal in $\Delta'(p,t)$.

  Suppose now that $\{f_n\}$ is not quasinormal at $p$.
  Then there exists points $z_j\in \Delta'(p,\delta)$ $(j=1,2,\ldots)$
  and a subsequence of $\{f_n\}$ (that we continue to call $\{f_n\}$)
  such that $z_j \rightarrow p$ as $j\rightarrow \infty$
  and no subsequence of $\{f_n\}$ is normal at any $z_j$, $j=1,2,\cdots$, see \cite[Thm. 4.4]{Shahara20014}.
  Let $E=\{z_j:j=1,2,\cdots\}$.
  Taking a subsequence of $\{f_n\}$ (that we continue to call $\{f_n\}$),
  we may assume that $f_n\overset{\chi}{\Longrightarrow} H$ in $\Delta'(p,\delta)\backslash E$.
  By Lemma  \ref{lquasichiefa}, we have $H'\equiv \psi$ and $H(z_j)=0$.
  It follows that $H$ is holomorphic in $\Delta'(p,t)$ and $H'\equiv \psi$ there.
  Moreover, since $\psi$ has no essential singularity at $p$, the same is true of $H$.
  But then $H(z_j)=0$ for $j=1,2,\ldots$ implies $H\equiv 0$, which contradict $H'\equiv \psi\not \equiv 0$.
 \end{proof}
\begin{lemma}\label{lnormalzz}
  Let $\{f_n\}$ be a family of meromorphic functions in $\Delta(z_0,r)$ and $k\neq 0$ an integer.
  Let $\{b_n\}$ be a sequence of holomorphic functions in $\Delta(z_0,r)$ such that $b_n\Rightarrow b$ in $\Delta(z_0,r)$,
  where $b(z)\neq 0$ is a holomorphic function in $\Delta(z_0,r)$.
   Suppose that $$f_n(z)\neq 0,\,f'_n(z)\neq (z-z_0)^kb_n(z),\quad z\in \Delta(z_0,r).$$
  Then $\{f_n\}$ is normal in $\Delta(z_0,r)$.
\end{lemma}
\begin{proof} Without loss of generality, we assume that $z_0=0$
and $r=1$.
   By Lemma \ref{guyongxing}, $\{f_n\}_{n=1}^\infty$ is normal in
   $\Delta' $. Without loss of generality,
   $f_n\overset\chi\Rightarrow f$ in $\Delta' ,$ and thus
   since $f_n\ne 0,$ then $\frac{1}{f_n}\Rightarrow \frac{1}{f}$
   in $\Delta' .$ The possibilities of the limit function
   $\frac{1}{f}$ are $\frac{1}{f}\equiv\infty$ or that
   $\frac{1}{f}$ is a holomorphic function. In the last case, we
   get by the maximum principle that
   $\{\frac{1}{f_n}\}_{n=1}^\infty$ converges to some holomorphic
   function in $\Delta $ and we are done. Hence we can
   assume that $\frac{1}{f}\equiv\infty,$ i.e., $f\equiv0$ in
   $\Delta' .$

   Obviously, we have
 \begin{equation}\nonumber
   f'_n(z)-z^kb_n(z)\Rightarrow -z^kb(z),\quad \text{in}
   \quad \Delta'.
 \end{equation}
 Since $f'_n(z)-z^kb_n(z)\neq 0$ in $\Delta$, we deduce that
 $$H_n(z)=\frac{1}{f'_n(z)-z^kb_n(z)}$$
 is holomorphic in $\Delta$, so by the maximum principle
 \begin{equation}\label{equationlemma3.10a}
   H_n(z)\Rightarrow -\frac{1}{z^kb(z)},\quad z\in \Delta.
 \end{equation}
 Obviously, $-\frac{1}{z^kb(z)}$ is an analytic function in $\Delta$ and hence $k$ is a negative integer.
 Let $m=-k$. By (\ref{equationlemma3.10a}) and Rouch\'{e}'s Theorem for holomorphic function, for sufficiently large $n$,
 \begin{equation}\label{equationlemma3.10b}
   n\left(\frac{1}{2}, f'_n(z)-\frac{b_n(z)}{z^m}\right)=n\left(\frac{1}{2}, \frac{1}{H_n(z)}\right)=
   n\left(\frac{1}{2}, -\frac{b(z)}{z^m}\right)=m.
 \end{equation}
 Since $f'_n(z)\neq \frac{b_n(z)}{z^m}$, then $f'_n(0) \neq \infty$ and hence $f_n(0) \neq \infty$.
 Obviously, $0$ is a pole of $f'_n(z)-\frac{b_n(z)}{z^m}$ of order $m$.
 By (\ref{equationlemma3.10b}), for sufficiently large $n$,
 $f_n$ has no poles in $\Delta(0,\frac{1}{2})$, i.e., $f_n$ is holomorphic in $\Delta(0,\frac{1}{2})$.
 We have proved that $f_n\Rightarrow 0$ in $\Delta'(z_0,r)$,
 so, by the maximum principle, $f_n\Rightarrow 0$ in $\Delta(0(\frac{1}{2}))$
 and  $\{f_n\}$ is   normal at $0$.
 \end{proof}
\begin{lemma}\label{lnormalpp}
  Let $\mathcal F=\{f_n\}$ be a family of holomorphic functions in $\Delta(z_0,r)$, all of whose zeros are multiple.
  Let $k \geq 2$ be a positive integer and $\{b_n\}$ be a sequence of holomorphic functions in $\Delta(z_0,r)$
  such that $b_n\Rightarrow b$ in $\Delta(z_0,r)$, where $b(z)\neq 0$ is a holomorphic function in $\Delta(z_0,r)$.
  Suppose that $$f'_n(z)\neq\frac{b_n(z)}{z^k},\quad z\in \Delta(z_0,r).$$
  Then $\mathcal F$ is normal in $\Delta(z_0,r)$.
\end{lemma}
\begin{proof} Without loss of generality, we assume that $z_0=0$
and $r=1$.
 By Lemma \ref{lnormala}, $\mathcal F$ is normal in $\Delta'$.
 Suppose that $\mathcal F$ is not normal at $0$.
 Let $\mathcal F_1=\{F_n\}$, where $F_n(z)=z^kf_n(z)$.
 Obviously, all zeros of $\mathcal F_1$ are multiple in $\Delta$.
 By Lemma \ref{lnonormalconvert}, $\mathcal F_1$ is not normal at $0$.
 By Lemma \ref{zp}, there exist points $z_n\rightarrow 0$, and positive numbers $\rho_n\rightarrow 0$
 and a subsequence of $\{F_n\}$ (that we continue to call $\{F_n\}$)
 such that
\begin{equation}\label{3.5a}g_n(\zeta)=\frac{F_n(z_n+\rho_n\zeta)}{\rho_n}=\frac{(z_n+\rho_n\zeta)^kf_n(z_n+\rho_n\zeta)}{\rho_n}\Rightarrow
g(\zeta),\quad \zeta\in\mathbb{C},\end{equation}
 where $g(\zeta)$ is a nonconstant holomorphic function in $\mathbb{C}$, all of whose zeros are multiple.

We consider the following two cases.

{\bf Case 1:} $z_n/\rho_n\rightarrow\infty$.

Observe first that since $\frac{z}{z_n+\rho_n\zeta}\Rightarrow 1$
in $\mathbb C,$ then \eqref{3.5a} is equivalent to
$\frac{z_n^k}{\rho_n}f_n(z_n+\rho_n\zeta)\Rightarrow g(\zeta)$ in
$\mathbb C.$ Differentiating gives
$z_n^kf_n'(z_n+\rho_n\zeta)\Rightarrow g'(\zeta)$ in $\mathbb C,$
and by the same reasoning as above it is equivalent to
$(z_n+\rho_n\zeta)^kf_n'(z_n+\rho_n\zeta)\Rightarrow g'(\zeta)$ in
$\mathbb C.$

Now we have
$$0\neq \frac{f'_n(z_n+\rho_n\zeta)}{\frac{b_n(z_n+\rho_n\zeta)}{(z_n+\rho_n\zeta)^k}}-1=\frac{(z_n+\rho_n\zeta)^k
f'_n(z_n+\rho_n\zeta)}{b_n(z_n+\rho_n\zeta)}-1 \Rightarrow
\frac{g'(\zeta)}{b(0)}-1\quad \text{in}\quad  \mathbb{C}. $$By
Hurwitz's Theorem, either $g'(\zeta)-b(0)\equiv 0$ in
$\mathbb{C}$, or $g'(\zeta)-b(0)\neq 0$ in $\mathbb{C}$. If
$g'(\zeta)-b(0)\equiv 0$, then $g(\zeta)=b(0)\zeta+c$ which is
impossible since  all the zeros of $g(\zeta)$ must be multiple. If
$g'(\zeta)-b(0)\neq 0$, then, by Lemma \ref{lrationlwang},
$g(\zeta)=\frac{b(0)(\zeta-a)^2}{\zeta-b}$ which contradicts that
$g(\zeta)$ is an entire function.

{\bf Case 2:}
$z_n/\rho_n\rightarrow\alpha\quad\text{in}\quad\mathbb C$

We have
$$G_n(\zeta)=\frac{F_n(\rho_n\zeta)}{\rho_n}=\frac{F_n(z_n+\rho_n(\zeta-\frac{z_n}{\rho_n}))}{\rho_n}\Rightarrow
g(\zeta-\alpha)\quad\text{in}\quad\mathbb{C}.$$ Writing
$$G(\zeta)=g(\zeta-\alpha)\quad\text{and}\quad\psi_n(\zeta)=\rho_n^{k-1}f_n(\rho_n\zeta)=\frac{F_n(\rho_n\zeta)}{\rho_n}\frac{1}
{\zeta^k},$$ we have
$$\psi_n(\zeta)\Rightarrow\frac{G(\zeta)}{\zeta^k}=\psi(\zeta),\quad\text{in}\quad  \mathbb{C}\setminus\{0\}.$$
Since for all $n$, $\psi_n(\zeta)$ is a holomorphic function, we have, by the maximum principle,
$$\psi_n(\zeta)\Rightarrow\psi(\zeta)\quad\text{in}
\quad \mathbb{C},$$ where $\psi(\zeta)$ is an entire function.
Hence we have
$$0\neq \rho_n^k \bigg(f'_n(\rho_n\zeta)-\frac{b_n(\rho_n\zeta)}
{(\rho_n\zeta)^k}\bigg) =\psi'_n(\zeta)-\frac{
b_n(\rho_n\zeta)}{\zeta^k} {\Longrightarrow }
 \psi'(\zeta)-\frac{b(0)}{\zeta^k}\quad\text{in}\quad \mathbb{C}\setminus\{0\}.$$
By Hurwitz's Theorem, either
$\psi'(\zeta)\equiv\frac{b(0)}{\zeta^k}$ in $\mathbb
C\setminus\{0\}$ (and hence in $\mathbb C$), or $\psi'(\zeta)\neq
\frac{b(0)}{\zeta^k}$ in $\mathbb C\setminus\{0\}$ (and in fact in
$\mathbb C$, since $b(0)/\zeta^k=\infty$ at $\zeta=0)$. Since
$\psi$ is a holomorphic function, the first alternative obviously
cannot hold. Thus $\psi'(\zeta)\neq \frac{b(0)}{\zeta^k}$. It then
follows from Theorem \ref{trational} and Lemma \ref{lrationalb}
that $\psi(\zeta)\equiv c$. Since $g$ is not a constant function
and $g(\zeta-\alpha)=\zeta^k\psi(\zeta)$, we have $c\neq0$. Now we
have
\begin{equation}\label{equationppa}
  \psi_n(\zeta)=\rho_n^{k-1}f_n(\rho_n\zeta)\Rightarrow c,\quad \text{in}
\quad \mathbb{C}.
\end{equation}
$$$$

Now suppose that there exists $\delta>0$ such that for
sufficiently large $n$, $f_n(z)\neq0$ for $z\in \Delta(0,\delta)$.
Then by Lemma \ref{lnormalzz}, $\{f_n\}$ is normal at $0$. A
contradiction.

  Otherwise, taking a subsequence and renumbering if necessary,
 we may assume that there exist $z_n^\ast\rightarrow0$ such that $f_n(z_n^\ast)=0$.
 We may assume that $z_n^\ast$ is the zero of $f_n$ of smallest modulus.
 By (\ref{equationppa}), we have, for sufficiently large $n$, $f_n(0)\neq 0$ and hence $z_n^\ast\neq 0$.
 By Hurwitz's Theorem and (\ref{equationppa}), we have $\frac{z_n^\ast}{\rho_n}\rightarrow\infty$.

Let $G_n(\zeta)=(z_n^\ast)^{k-1}f_n(z_n^\ast\zeta)$. Then we have
$${G'_n}(\zeta)=(z_n^\ast)^kf'_n(z_n^\ast\zeta)\neq \frac{b(z_n^\ast\zeta)}{\zeta^k}.$$
Since $G_n(\zeta)\neq 0$ in $\Delta$, it follows from Lemma
\ref{lnormalzz} that $\{G_n\}$ is normal in $\Delta$. By Lemma
\ref{lquasichiefb}, $\{G_n\}$ is quasinormal in $\mathbb{C}$.
Thus, there exists a subsequence of $\{G_n \}$ (that we continue
to call $\{G_n \}$) and $E\subset \mathbb{C}$ such that
\begin{lista}
\item [$(b1)$] $E$ have no accumulation point in $\mathbb{C}$;
\item [$(b2)$] $G_n(\zeta)\overset{\chi}{\Longrightarrow}G(\zeta)$
in $\mathbb{C}\backslash E$; \item [$(b3)$] for each $\zeta_0\in
E$, no subsequence of $\{G_n \}$ is normal at $\zeta_0$.
\end{lista}
Obviously, we have $E\cap\Delta=\emptyset$ and all zeros of $G(\zeta)$ are multiple.

Since $G_n(0)=\left(\frac{z_n^\ast}{\rho_n}
\right)^{k-1}\psi_n(0)\rightarrow \infty$ and $\{G_n(\zeta)\}$ is
a family of holomorphic functions, we have $G(\zeta)\equiv \infty$
in $C\backslash E$. Suppose that  $1\not\in E$. Since $G_n(1)=0$,
we have $G(1)=0$ which contradicts that $G(\zeta)\equiv \infty$ on
$C\backslash E$. Thus we have $1\in E$. By Lemma
\ref{lquasichiefa}, $G(\zeta)=\int_1^\zeta\frac{b(0)}{\xi^k}d\xi$
on $C\backslash E$, which also contradicts that $G(\zeta)\equiv
\infty$ in $\mathbb{C}\backslash E$.
\end{proof}
\section{Proof of Theorem}
\begin{proof}
We assume that $f'=h$ has at most finitely many zeros and derive a contradiction.

We claim that there exist $t_n\rightarrow\infty$ and $\varepsilon_n \rightarrow 0$ such that
\begin{equation}\label{sphericalcover}
  S(\Delta(t_n,\varepsilon_n),f)=\frac{1}{\pi}\iint_{|z-t_n|<\varepsilon_n}[f^\#(z)]^2dxdy\rightarrow \infty.
\end{equation}
Otherwise there would exist $\varepsilon >0$ and $M >0 $ such that for all $z_0\in \mathbb{C}$, we have
$$S(\Delta(z_0,\varepsilon),f)=\frac{1}{\pi}\iint_{|z-z_0|<\varepsilon}[f^\#(z)]^2dxdy < M.$$
Then $S(r,f)=\frac{1}{\pi}\iint_{|z|<r}[f^\#(z)]^2dxdy=O(r^2)$.
Thus
$$T_0(r,f)=\int_0^r\frac{S(t)}{t}=O(r^2).$$
On the other hand, by Lemma \ref{lelliptic},
$T_0(r,h)=Ar^2(1+(1))$ as  $r\rightarrow \infty$ where $A>0$ is a
constant. Hence it is impossible that $T_0(r,h)=o\{T_0(r,f)\}$ and
\eqref{sphericalcover} follows.

Let $\omega_1,\omega_2$ be the two fundamental periods of $h(z)$
and $P(0 \in P)$ be a fundamental parallelogram of $h(z)$. There
exist integers $i_n$ and $j_n$ such that $z_n\in P$, where
$z_n=t_n-i_n\omega_1-j_n\omega_2$. There exists a subsequence of
$\{z_n\}$ (that we continue to call $\{z_n\}$) such that
$z_n\rightarrow z_0$ as $n\rightarrow \infty$. Let
$f_n(z)=f(z+i_n\omega_1+j_n\omega_2)$. By (\ref{sphericalcover})
we have
\begin{equation}\label{sphericalcovera}
  S(\Delta(z_n,\varepsilon_n),f_n)=S(\Delta(t_n,\varepsilon_n),f)\rightarrow \infty,
\end{equation}
and hence, there exists  $z_n^\ast$ ($z_n^\ast\rightarrow z_0$)
such that $f_n^\#(z_n^\ast)\rightarrow \infty$ as $n\rightarrow \infty$.
Without loss of generality, we can assume that $z_0=0$.
Hence we have that no subsequence of $\{f_n\}$ is normal at $0$ and by (\ref{sphericalcovera}),
\begin{equation}\label{sphericalcoverb}
  S(\Delta(z_n,\varepsilon_n),f_n)\rightarrow \infty,
\end{equation}
where $z_n\rightarrow 0$ and $\varepsilon_n\rightarrow 0$ as $n\rightarrow\infty$.

 There exists $R>0$ such that $\overline P\subset \Delta(0,R)$ and $\Delta(z_n,\varepsilon_n)\subset \Delta(0,R)$ for each $n$. Set $D=\Delta(0,R)$. Obviously, we have $z_0\in D$.
 By assumption, for sufficiently large $n$,
 $$f'_n(z)=f'(z+i_n\omega_1+j_n\omega_2)\neq h(z+i_n\omega_1+j_n\omega_2)=h(z),\quad z\in D.$$
 Without loss of generality, we can assume that for all $n\in \mathbb{N}$,
 $f'_n(z)\neq h(z)$ in $D$.

 Now, $\{f_n\}$ is a family of functions meromorphic in $D$ such that
\begin{lista}
\item [$(a^\ast)$] all zeros of $\{f_n\}$ are multiple in $D$;
\item [$(b^\ast)$] for each $n$, $f'_n(z)\neq h(z)$ in $D$, where
$h$ is a nonconstant elliptic function; \item [$(c^\ast)$] no
subsequence of $\{f_n\}$ is normal at $0$.
\end{lista}
It follows from Lemma \ref{lquasichiefb} that $\{f_n\}$ is
quasinormal in $D$. Hence there exists $\tau>0$ such that
$\{f_n\}$ is normal in $\Delta'(0,\tau)$ and $h(z)\neq 0$ on
$\Delta'(0,\tau)$.

Without loss of generality, we may assume that $\tau=1$. Then
there exists a subsequence of $\{f_n\}$ (that we continue to call
$\{f_n\}$) such that
\begin{lista}
\item [$(a)$] all zeros of $\{f_n\}$ are multiple in $\Delta$;
\item [$(b)$] for each $n$, $f'_n(z)\neq h(z)$ in $\Delta$, where
$h\neq 0$ in $\Delta'$; \item [$(c)$] no subsequence of $\{f_n\}$
is normal at $0$. \item [$(d)$] $f_n(z)\Rightarrow f_0(z)$ on
$\Delta'(0,1)$
\end{lista}

{\bf Case 1:} $h(0)\neq0,\infty$.

By Lemma \ref{lquasichiefa}, $f_0(z)=\int_{\zeta=0}^z
h(\zeta)d\zeta$ and for sufficiently large $n$, there exists
$1>\delta>0$ such that $f_n$ has a single zero of order $2$ in
$\Delta(0,\delta)$. By Lemma \ref{lalforscoverb}, there exists
$M>0$ such that $S(\frac{\delta}{2},f_n)<M$ which contracts (\ref
{sphericalcoverb}).

{\bf Case 2:} $h(0)=0$.

 Suppose that $0$ is a zero of order $k$ of $h(z)$, where $k$ is a positive integer.
 Let us assume, making standard normalizations, that for $z\in \Delta$
 $$h(z)=z^k+a_{k+1}z^{k+1}+\cdots=z^k\widehat h(z)$$
 where $\widehat h(z)\neq0,\infty$ in $\Delta(0,1)$ and $\widehat h(0)=1$.

 We claim that for each $\delta>0$,
 there exists at least one zero of $f_n$ in $\Delta(0,\delta)$  for sufficiently large $n$.
 Otherwise, there exists a subsequence of $\{f_n\}$ (that we continue to call $\{f_n\}$) such that $f_n(z)\neq 0$ in $\Delta(0,\delta)$. By Lemma \ref{lschwick}, $\{f_n\}$ is normal at $0$.
 A contradiction.

Hence, taking a subsequence and renumbering if necessary,
we may assume that $a_n\rightarrow 0$ is the zero of $f_n$ of smallest modulus.
Since $f'_n(z)\neq h(z)$ and all the zeros of $\{f_n\}$ are multiple,
we have $f_n(0)\neq 0$ and hence $a_n\neq 0$.
Let $F_n=\frac{f_n(a_n\zeta)}{a_n^{k+1}}$. We have that
\begin{lista}
\item [$(a1)$] $F_n(\zeta)\neq 0$ in $\Delta$; \item [$(a2)$] all
zeros of $F_n(\zeta)$ are multiple; \item [$(a3)$]
$F'_n(\zeta)\neq \zeta^k\widehat h(a_n\zeta)$  and $F_n(1)=0$.
\end{lista}
By Lemma \ref{lnormalzz}, $\{F_n(\zeta)\}$ is normal in $\Delta$.
By Lemma \ref{lquasichiefb}, $\{F_n(\zeta)\}$ is quasinormal on
$\mathbb{C}$. Thus, there exists a subsequence of $\{F_n(\zeta)\}$
(that we continue to call $\{F_n(\zeta)\}$) and $E_1 \subset
\mathbb{C}$ such that
\begin{lista}
\item [$(b1)$] $E_1$ have no accumulation point in $\mathbb{C}$;
\item [$(b2)$] $F_n(\zeta)\overset{\chi}{\Longrightarrow}F(\zeta)$
on $\mathbb{C}\backslash E_1$; \item [$(b3)$] for each $\zeta_0\in
E_1$, no subsequence of $\{F_n(\zeta)\}$ is normal at $\zeta_0$.
\end{lista}
Obviously, $E_1\cap\bigtriangleup=\emptyset$ and all zeros of
$F(\zeta)$ are multiple in $\mathbb{C}\backslash E_1$.

{\bf Case 2.1:} $1\not\in E_1$

Since all zeros of $\{F_n(\zeta)\}$ are multiple and $F_n(1)=0$,
we have that $F(1)=F'(1)=0$ (recall that $\hat h(0)=1$) and hence
$F(\zeta)$ is a meromorpic function in $\mathbb{C}\backslash E_1$.

We claim that $E_1=\emptyset$. Otherwise, let $\zeta_0\in E_1$,
Obviously, $\zeta_0\neq 0$. By Lemma \ref{lquasichiefa},
$F'(\zeta)=\zeta^k$. Recall that $\hat h(0)=1$ and hence $F'(1)=1$
which contradicts that $F'(1)=0$.

By Lemma \ref{lhurwitz}, either $F'(\zeta)\equiv \zeta^k$ in
$\mathbb{C}$, or $F'(\zeta)\neq\zeta^k$ in $\mathbb{C}$. If
$F'(\zeta)\equiv \zeta^k$ in $\mathbb{C}$, then $F'(1)=1$ which
contradicts that $F'(1)=0$. If $F'(\zeta)\neq\zeta^k$ on
$\mathbb{C}$. By Theorem \ref{trational}, $F$ must be rational and
then by Lemma \ref{lrationalc},
$$F(\zeta)=\frac{\prod\limits_{i=1}^{m+k+1}(\zeta-\alpha_i)}{(k+1)(\zeta-\beta)^m},$$
 where    $m$ is a nonnegative integer, $\beta\in \mathbb{C}$, $\alpha_i\neq0,\beta$, $1\le i\le m+k+1$.

Hence, we have
\begin{equation}\label{equation21a}
  F_n(\zeta)\overset{\chi}{\Longrightarrow}\frac{\prod\limits_{i=1}^{m+k+1}(\zeta-\alpha_i)}{(k+1)(\zeta-\beta)^m},\quad{in}\quad \mathbb{C}
\end{equation}

By Hurwitz's Theorem, there exist sequences
$\zeta_{n,i}\rightarrow \alpha_i$ and $\eta_{n,j}\rightarrow
\beta$  as $n\to\infty$ (counting multiplicities of zeros and
poles, respectively), such that for sufficiently large $n$,
$F_n(\zeta_{n,i})=0$ and $F_n(\eta_{n,j})=\infty$, where
$i=1,2,\cdots, m+k+1$ and $j=1,2,\cdots, m$. Writing
$z_{n,i}=a_n\zeta_{n,i}$. Thus, $f_n(z_{n,i})=0$ and
$z_{n,i}\rightarrow 0$ as $n\rightarrow\infty$, where
$i=1,2,\cdots, m+k+1$. Set
$B_n=\{z_{n,1},z_{n,2},\cdots,z_{n,m+k+1}\}$.

 We claim that for each $\delta>0$,
 there exists at least $m+k+2$ zeros of $f_n$ in $\Delta(0,\delta)$  for sufficiently large $n$.
 Otherwise, there exists a subsequence of $\{f_n\}$ (that we continue to call $\{f_n\}$) such that $f_n(z)$ have $m+k+1$  zeros in $\Delta(0,\delta)$.
 If $f_0(z)\equiv \infty$, then by Lemma \ref{lalforscoverb}, there exists $M>0$ such that for $n$ sufficiently large, $S(\frac{\delta}{2}, f_n)<M$ which contradicts (\ref{sphericalcoverb}).
 If $f_0(z)\not\equiv \infty$, then by Lemma \ref{ldnormal},
 $f_0(z)=\frac{\zeta^{k+1}+c}{k+1}$ in $\Delta'$, where $c$ is a constant. By Lemma \ref{lalforscoverb}, there exists $M>0$ such that for $n$ sufficiently large, $S(\frac{\delta}{2}, f_n)<M$ which also contradicts (\ref{sphericalcoverb}).

 Hence, taking a subsequence and renumbering if necessary,
 we may assume that $b_n\rightarrow 0$ is the zero of $f_n$ of smallest modulus in $\Delta\backslash B_n$.
 Set $r_n=\frac{a_n}{b_n}$. Then we have $F_n(\frac{1}{r_n})=0$.
 Since $b_n\not\in B_n$, Obviously $\frac{1}{r_n}\neq\zeta_{n,i}$, where $i=1,2,\cdots, m+k+1$.
 By Hurwitz's Theorem and (\ref{equation21a}),
 we have $\frac{1}{r_n}\rightarrow\infty$ and hence $r_n\rightarrow 0$ as $n\rightarrow\infty$.

Let $G_n(\zeta)=\frac{f_n(b_n\zeta)}{b_n^{k+1}}$. We have that for sufficiently large $n$,
\begin{lista}
\item [$(d1)$] $G_n(\zeta)$ have only $m+k+1$ zeros
$r_n\zeta_{n,i}$ in $\Delta$. Obviously,
$|r_n\zeta_{n,i}|\rightarrow 0$, as $n\rightarrow\infty$, where
$i=1,2,\cdots, m+k+1$; \item [$(d2)$] all zeros of $G_n(\zeta)$
are multiple; \item [$(d3)$] $G'_n(\zeta)\neq \zeta^k\widehat
h(b_n\zeta)$ and $G_n(1)=0$.
\end{lista}
By Lemma \ref{lnormalzz}, $\{G_n(\zeta)\}$ is normal in $\Delta'$.
By Lemma \ref{lquasichiefb}, $\{G_n(\zeta)\}$ is quasinormal on
$\mathbb{C}$. Thus, there exists a subsequence of $\{G_n(\zeta)\}$
(that we continue to call $\{G_n(\zeta)\}$) and $E_2 \subset
\mathbb{C}$ such that
\begin{lista}
\item [$(e1)$] $E_2$ have no accumulation point in $\mathbb{C}$;
\item [$(e2)$] $G_n(\zeta)\overset{\chi}{\Longrightarrow}G(\zeta)$
on $\mathbb{C}\backslash E_2$; \item [$(e3)$] for each $\zeta_0\in
E_2$, no subsequence of $\{G_n(\zeta)\}$ is normal at $\zeta_0$.
\end{lista}
Obviously, $E_2\cap\Delta'=\emptyset$ and all zeros of $G(\zeta)$
are multiple in $\mathbb{C}\backslash E_2$.

Let $$G_n^\ast(\zeta)=G_n(\zeta)\frac{\prod\limits_{j=1}^{m}(\zeta-r_n\eta_{n,j})}{\prod\limits_{i=1}^{m+k+1}(\zeta-r_n\zeta_{n,i})},$$
$$F_n^\ast(\zeta)=F_n(\zeta)\frac{\prod\limits_{j=1}^{m}(\zeta-\eta_{n,j})}{\prod\limits_{i=1}^{m+k+1}(\zeta-\zeta_{n,i})}.$$
By (\ref{equation21a}),
\begin{equation} \nonumber
 G_n^\ast(r_n\zeta)=F_n^\ast(\zeta)\Rightarrow \frac{1}{k+1},\quad \zeta\in \mathbb{C}
\end{equation}
Hence
\begin{equation}\label{equation21b}
G_n^\ast(0)\rightarrow\frac{1}{k+1}.
\end{equation}
We claim that $G(\zeta)$ is a meromorphic function on
$\mathbb{C}\backslash E_2$. Suppose that $G(\zeta)\equiv \infty$.
Obviously, $G_n^\ast(\zeta)$ have no zeros in $\Delta$. Applying
the maximum principle to the sequence $\frac{1}{G_n^\ast(\zeta)}$
of analytic functions, we see that $G_n^\ast(\zeta)\equiv \infty$
in $\Delta $ which contradict $G_n^\ast(0)\Rightarrow
\frac{1}{k+1}$.

We claim that $G(\zeta)=\frac{\zeta^{k+1}}{k+1}$ on
$\mathbb{C}\backslash E_2$. Indeed, since $G(\zeta)$ is a
meromorphic function in $\mathbb{C}\backslash E_2$, then by Lemma
\ref{ldnormal}, $G(\zeta)=\frac{\zeta^{k+1}+c}{k+1}$, where $c$ is
a constant. Since $G_n^\ast(\zeta)$ have no zeros in $\Delta$, we
have
$$G_n^\ast(\zeta)\Rightarrow \frac{\zeta^{k+1}+c}{k+1}\frac{1}{\zeta^{k+1}},\quad \zeta\in \Delta.$$
Hence, $G_n^\ast(0)\rightarrow
\left.\frac{\zeta^{k+1}+c}{k+1}\frac{1}{\zeta^{k+1}}\right|_{\zeta=0}$
. By (\ref{equation21b}), we get that $c=0$.

Suppose that $1\not\in E_2$. Since $G_n(1)=0$, we have $G(1)=0$
which contradicts that $G(\zeta)=\frac{\zeta^{k+1}}{k+1}$. Thus,
$1\in E_2$. By Lemma \ref{lquasichiefa},
$G(\zeta)=\int_1^\zeta\xi^{k}d\xi=\frac{\zeta^{k+1}-1}{k+1}$ which
contradicts that $G(\zeta)=\frac{\zeta^{k+1}}{k+1}$.

{\bf Case 2.2:} $1\in E_1$

By Lemma \ref{lquasichiefa},
\begin{equation}\label{equation22a}
 F(\zeta)=\int_1^{\zeta}\xi^k d\xi=\frac{\zeta^{k+1}-1}{k+1},\quad \zeta\in \mathbb{C}\backslash E_1.
\end{equation}
Let $e_j$ be the $j$th root of the equation $\zeta^{k+1}-1=0$,
$j=1,2,\cdots k+1$.

\begin{claim*}
  $E_1=\{e_1,e_2,\cdots,e_{k+1}\}$.
\end{claim*}

\begin{proof}
 Suppose that $\zeta_0\not\in E_1$, where $\zeta_0^{k+1}-1=0$.
Obviously, $\zeta_0$ is a zero of $F(\zeta)$ but not a multiple
zero of $F(\zeta)$ which contradicts that all of zeros of
$F(\zeta)$ are multiple.

Suppose that $\zeta_0\in E_1$, where $\zeta_0^{k+1}-1\neq 0$.
 By Lemma \ref{lquasichiefa}, $F(\zeta)=\int_{\zeta=\zeta_0}^{\zeta}\xi^{k+1} d\xi=\frac{\zeta^{k+1}-\zeta_0^{k+1}}{k+1}$,
 $\zeta\in\mathbb{C}\backslash E_1$. By (\ref{equation22a}), $\zeta_0^{k+1}=1$, a contradiction.
 \end{proof}

By Lemma \ref{lquasichiefa}, there exists $\delta_j>0$ such that
for sufficiently large $n$, $F_n(\zeta)$ have a single zero
$\zeta_{n,j}\rightarrow e_j$ of order 2 and a single pole
$\eta_{n,j}\rightarrow e_j$ of order 1 in $\Delta(e_j,\delta_j)$.

Set $z_{n,j}=a_n\zeta_{n,j}$. Thus, $f_n(z_{n,j})=0$ and
$z_{n,j}\rightarrow 0$ as $n\rightarrow\infty$, where
$j=1,2,\cdots, k+1$. Set
$B_n=\{z_{n,1},z_{n,2},\cdots,z_{n,k+1}\}$.
 Do as in   Case 2.1,
 we may assume that $b_n\rightarrow 0$ is the zero of $f_n$ of smallest modulus in $\Delta\backslash B_n$.
 Set $r_n=\frac{a_n}{b_n}$. Obviously, $F_n(\frac{1}{r_n})=0$.
 Since $b_n\not \in B_n$, $\frac{1}{r_n}\neq \zeta_{n,j}$, where $j=1,2,\cdots, k+1$.
 Since $F_n(\zeta)$ have a single zero $\zeta_{n,j}\rightarrow e_j$ of order 2 in $\Delta(e_j,\delta_j)$,
 by Hurwitz's Theorem and (\ref{equation22a})
 we have $\frac{1}{r_n}\rightarrow\infty$ and hence $r_n\rightarrow 0$ as $n\rightarrow\infty$.

Let $G_n(\zeta)=\frac{f_n(b_n\zeta)}{b_n^{k+1}}$. We have that for sufficiently large $n$,
\begin{lista}
\item [$(f1)$] $G_n(\zeta)$ have only $k+1$ zeros $r_n\zeta_{n,i}$
of order 2 and at least $k+1$ poles $r_n\eta_{n,i}$ of order 1 in
$\Delta$. Obviously,  $|r_n\zeta_{n,i}|\rightarrow 0$ and
$|r_n\eta_{n,i}|\rightarrow 0$ as $n\rightarrow\infty$; \item
[$(f2)$] all zeros of $G_n(\zeta)$ are multiple; \item [$(f3)$]
$G'_n(\zeta)\neq \zeta^k\widehat h(b_n\zeta)$ and $G_n(1)=0$.
\end{lista}
By Lemma \ref{guyongxing}, $\{G_n \}$ is normal in $\Delta'$. By
Lemma \ref{lquasichiefb}, $\{G_n \}$ is quasinormal in
$\mathbb{C}$. Thus, there exist a subsequence of $\{G_n(\zeta)\}$
(that we continue to call $\{G_n(\zeta)\}$) and $E_3 \subset
\mathbb{C}$ such that
\begin{lista}
\item [$(g1)$] $E_3$ have no accumulation point in $\mathbb{C}$;
\item [$(g2)$] $G_n(\zeta)\overset{\chi}{\Longrightarrow}G(\zeta)$
on $\mathbb{C}\backslash E_3$; \item [$(g3)$] for each $\zeta_0\in
E_3$, no subsequence of $\{G_n(\zeta)\}$ is normal at $\zeta_0$.
\end{lista}
Obviously, $E_3\cap\Delta'=\emptyset$  and  all zeros of
$G(\zeta)$ are multiple in $\mathbb{C}\backslash E_3$;

Let
$$G_n^\ast(\zeta)=G_n(\zeta)\frac{\prod\limits_{j=1}^{k+1}(\zeta-r_n\eta_{n,j})}{\prod\limits_{j=1}^{k+1}(\zeta-r_n\zeta_{n,j})^2},$$
$$F_n^\ast(\zeta)=F_n(\zeta)\frac{\prod\limits_{j=1}^{k+1}(\zeta-\eta_{n,j})}{\prod\limits_{j=1}^{k+1}(\zeta-\zeta_{n,j})^2}.$$

By (\ref{equation22a}),
\begin{equation} \nonumber
 G_n^\ast(r_n\zeta)=F_n^\ast(\zeta)\Rightarrow \frac{1}{k+1},\quad \zeta\in \mathbb{C}.
\end{equation}
Hence
\begin{equation}\label{equation22b}
  G_n^\ast(0)\rightarrow\frac{1}{k+1}.
\end{equation}
We claim that $G(\zeta)$ is a meromorphic function on
$\mathbb{C}\backslash E_3$. Suppose that $G(\zeta)\equiv \infty$.
Obviously, $G_n^\ast(\zeta)$ have no zeros in $\Delta$. Applying
the maximum principle to the sequence $\frac{1}{G_n^\ast(\zeta)}$
of analytic functions, we see that $G_n^\ast(\zeta)\equiv \infty$
in $\Delta $ which contradict $G_n^\ast(0)\Rightarrow
\frac{1}{k+1}$.

We claim that $G(\zeta)=\frac{\zeta^{k+1}}{k+1}$ on
$\mathbb{C}\backslash E_3$. Since $G(\zeta)$ is a meromorphic
function in $\mathbb{C}\backslash E_3$, by Lemma \ref{ldnormal},
$G(\zeta)=\frac{\zeta^{k+1}+c}{k+1}$, where $c$ is a constant.
Since $G_n^\ast(\zeta)$ have no zeros in $\Delta$, we have
$$G_n^\ast(\zeta)\Rightarrow \frac{\zeta^{k+1}+c}{k+1}\frac{1}{\zeta^{k+1}},\quad \zeta\in \Delta$$
If $c\ne0$ then we get  $G_n^\ast(0)\rightarrow\infty,$ and this
contradicts  \eqref{equation22b}.

We claim that $1\in E_3.$ Indeed, suppose that $1\not\in E_3$.
Since $G_n(1)=0$, we have $G(1)=0$ which contradicts that
$G(\zeta)=\frac{\zeta^{k+1}}{k+1}$. Thus, $1\in E_3$. By Lemma
\ref{lquasichiefa}, $G(\zeta)=\int_1^\zeta\zeta^{k
}d\xi=\frac{\zeta^{k+1}-1}{k+1}$ which contradicts that
$G(\zeta)=\frac{\zeta^{k+1}}{k+1}$.

{\bf Case 3:} $h(0)=\infty$.

  Suppose that $0$ is a pole of order $k$ of $h(z)$.
  Let us assume, making standard normalizations,
  that for $z\in \Delta$
  $$h(z)=\frac{1}{z^k}+\frac{a_{k-1}}{z^{k-1}}+\cdots=\frac{\widehat h(z)}{z^k},$$
  where $\widehat h(z)\neq0,\infty$ in $\Delta$ and $\widehat h(0)=1$.

  {\bf Case 3.1:} $k=1$.

  By Lemma \ref{zp} (with $\alpha=0$),
  there exists points $z_{n}\rightarrow 0$, positive numbers $\rho_{n}\rightarrow 0$
  and a subsequence of $\{f_{n}\}$(that we continue to call $\{f_{n}\}$) such that
  \begin{equation}\label{equation310}
    g_n(\zeta)=f_n(z_n+\rho_n\zeta)\overset{\chi}{\Longrightarrow} g(\zeta), \quad \zeta\in\mathbb{C},
  \end{equation}
    where $g(\zeta)$ is a nonconstant meromorphic function in $\mathbb{C}$, all of whose zeros are multiple.

Again we consider two cases.

 {\bf Case 1:}    $z_n/\rho_n\rightarrow\infty$.

 Consider
\begin{equation}\label{equation311a}
  \varphi_n(\zeta)=f_n(z_n+z_n\zeta)=f_n(z_n(1+\zeta)).
\end{equation}
  $$$$
  Since $g$ is nonconstant, there exist $\zeta_1$, $\zeta_2\in \mathbb{C}$ such that $g(\zeta_1)\neq g(\zeta_2)$.
We have for $j=1,2$,
$$g(\zeta_j)=\lim_{n\rightarrow\infty}f_n(z_n+\rho_n\zeta_j)=\lim_{n\rightarrow\infty}f_n(z_n+z_n(\frac{\rho_n}{z_n}\zeta_j))=\lim_{n\rightarrow\infty}\varphi_n(\frac{\rho_n}{z_n}\zeta_j).$$
Since $(\frac{\rho_n}{z_n})\zeta_j\rightarrow 0$ as $n\rightarrow\infty$,
the family $\{\varphi_n\}$ is not equicontinuous at $0$
and hence no subsequence of $\{\varphi_n\}$ is normal at $0$.
By (\ref{equation311a}),
\begin{equation}\label{equation31a}
  \varphi_n'(\zeta)=z_nf_n'(z_n(1+\zeta))\neq \frac{\widehat h(z_n(1+\zeta))}{1+\zeta}.
\end{equation}
Since all zeros of $\varphi_n$ are multiple and by
(\ref{equation31a}), the family $\{\varphi_n\}$ is quasinormal in
$\mathbb{C}$ by Lemma \ref {lquasichiefb}. Hence there exist a set
$E$ and a subsequence of $\{\varphi_n\}$ (that we continue to
denote by $\{\varphi_n\}$) such that
\begin{lista}
 \item [$(h1)$] $E$ has no accumulation in $\mathbb{C}$;
 \item [$(h2)$] $\varphi_n(\zeta)\overset{\chi}{\Longrightarrow}\varphi(\zeta)$ in $\mathbb{C}\backslash E$
                and all the zeros of $\varphi_n(\zeta)$ are multiple;
 \item [$(h3)$] for each $\zeta_0\in E$, no subsequence of $\{\varphi_n\}$ is normal at $\zeta_0$.
\end{lista}

 Obviously, $0\in E$.
 by Lemma \ref{lquasichiefa}, $\varphi(\zeta)=\int^\zeta_{0}\frac{1}{1+\xi}d\xi$. Hence, $\int^\zeta_{0}\frac{1}{1+\xi}d\xi$ is a multi-valued function in $\mathbb{C}\backslash E$.
 A contradiction.

 {\bf Case 2:}    $\frac{z_n}{\rho_n}\rightarrow\alpha \in\mathbb
 C.$

    Since $g_n'(\zeta)\neq \frac{\rho_n\widehat h(z_n+\rho_n\zeta)}{z_n+\rho_n\zeta}$ and $\frac{\rho_n\widehat h(z_n+\rho_n\zeta)}{z_n+\rho_n\zeta}\overset{\chi}{\Longrightarrow}\frac{1}{\zeta+\alpha}$ in
    $\mathbb{C}$, then
    by Lemma \ref{lhurwitz}, either $g'(\zeta)\equiv\frac{1}{\alpha+\zeta}$ in $\mathbb{C}$, or $g'(\zeta)\neq\frac{1}{\alpha+\zeta}$ in $\mathbb{C}$.
    However, since all poles of $g'$ are multiple, the first alternative obviously cannot hold,
    Thus  $g'(\zeta)\neq\frac{1}{\alpha+\zeta}$. By Theorem \ref{trational} we deduce that $g$ is rational.
    By the fundamental theorem of algebra, $g$ cannot be a polynomial.
    hence
    $$g'(\zeta)=\frac{1}{a+\zeta}+\frac{1}{p(\zeta)},$$
    where $p(\zeta)$ is a polynomial. It is easy to see that if
    $\deg p\le 1,$ then the right hand side  above cannot be
    a derivative of a nonconstant meromorphic function in $\mathbb
    C.$
    But then $$g'(\zeta)=\frac{1}{\zeta}+O\left(\frac{1}{|\zeta|^2}\right),\quad \zeta\to\infty,$$
    so that
    $$\frac{1}{2\pi i}\int_{\Gamma(0,R)}g'(\zeta)d\zeta=1+o(  R),\quad  R\to\infty$$
    for $R$ sufficiently large,
    where $\Gamma(0,R)$ is the positive oriented circle of radius $R$ about the origin.
    On the other hand, the  above integral must vanish, as it is the integral of a derivative over a closed curve. A contradiction.

{\bf Case 3.2:} $k\geq2$.

 We claim that for each $\delta>0$,
 there exists at least one pole of $f_n$ in $\Delta(0,\delta)$  for sufficiently large $n$.
 Otherwise, there exists a subsequence of $\{f_n\}$ (that we continue to call $\{f_n\}$)
 such that $\{f_n(z)\}$ is a sequence of functions holomorphic in $\Delta(0,\delta)$.
 By Lemma \ref{lnormalpp}, $\{f_n\}$ is normal at $0$. A contradiction.

 Hence, taking a subsequence and renumbering if necessary,
 we may assume that $a_n\rightarrow 0$ is the pole of $f_n(z)$ of smallest modulus.
 Since $f'_n(z)\neq h(z)$ and $h(0)=\infty$, we have $f(0)\neq \infty$ and hence $a_n\neq 0$.
 Let $F_n(\zeta)=a_n^{k-1}f_n(a_n\zeta)$, we have
\begin{lista}
\item [$(i1)$] $F_n(\zeta)$ is holomorphic function in $\Delta$;
\item [$(i2)$] all zeros of $F_n(\zeta)$ are multiple; \item
[$(i3)$] $F'_n(\zeta)\neq \frac{\widehat h(a_n\zeta)}{\zeta^k}$
and $F_n(1)=\infty$.
\end{lista}
By Lemma \ref{lnormalpp}, $\{F_n(\zeta)\}$ is normal at $\Delta$.
By Lemma \ref{lquasichiefb}, $\{F_n(\zeta)\}$ is quasinormal in
$\mathbb{C}$. Thus, there exists a subsequence of $\{F_n(\zeta)\}$
(that we continue to call $\{F_n(\zeta)\}$) and $E_4 \subset
\mathbb{C}$ such that
\begin{lista}
\item [$(j1)$] $E$ have no accumulation point in $\mathbb{C}$;
\item [$(j2)$] $F_n(\zeta)\overset{\chi}{\Longrightarrow}F(\zeta)$
in $\mathbb{C}\backslash E_4$; \item [$(j3)$]  for each
$\zeta_0\in E_4$, no subsequence of $\{F_n(\zeta)\}$ is normal at
$\zeta_0$.
\end{lista}
Obviously, $E_4\cap\Delta=\emptyset$ and all zeros of $F(\zeta)$
are multiple in $\mathbb{C}\backslash E_4$.

\begin{claim*} For $n$ sufficiently large,  $F_n(0)\neq 0$
(and hence $f_n(0)\neq 0$).
\end{claim*}
\begin{proof}
Otherwise, there exists a subsequence of $\{F_n(\zeta)\}$ (that we
continue to call $\{F_n(\zeta)\}$) such that $F_n(0)=0$. Thus
$F(0)=0$ and hence $F(\zeta)$ is a meromorphic function in
$\mathbb{C}\backslash E_4$. Suppose that $E_4$ is the empty set.
Since $F'_n(\zeta)\neq \frac{\widehat h(a_n\zeta)}{\zeta^k}$ and
$\frac{\widehat
h(a_n\zeta)}{\zeta^k}\overset{\chi}{\Longrightarrow}\frac{1}{\zeta^k}$
in $\mathbb{C}$, by Lemma \ref{lhurwitz}  we have either
$F'(\zeta)\equiv  \frac{1}{\zeta^k}$ in $\mathbb{C}$, or
$F'(\zeta)\neq\frac{1}{\zeta^k}$ in $\mathbb{C}$. If
$F'(\zeta)-\frac{1}{\zeta^k}\equiv
 0$ in $\mathbb{C}$, then we have $F(0)=\infty$ which contradicts that $F(0)=0$.
If $F'(\zeta)\neq\frac{1}{\zeta^k}$ in $\mathbb{C}$, then, by
Theorem \ref{trational} and  Lemma \ref{lrationalb}, we have
$F(\zeta)=c$, where $c\in \mathbb{C}$ is a constant. On the other
hand, since $F_n(1)=\infty$, we have $F(1)=\infty$, a
contradiction.

Thus, we can assume that $E_4$ is not the empty set. Suppose that
$\zeta_0\in E_4$. Obviously, $\zeta_0\neq 0$. By Lemma
\ref{lquasichiefa}, we have
$F(\zeta)=\int_{\zeta_0}^\zeta\frac{1}{\xi^k}d\xi$ which implies
$F(0)=\infty$, a contradiction to the fact that
$F_n(0)=0$.\end{proof}

Without loss of generality, we assume that  for each
$n\in\mathbb{N}$, $f_n(0)\neq 0$.

 We claim that for each $\delta>0$,
 there exists at least one zero of $f_n$ in $\Delta(0,\delta)$  for sufficiently large $n$.
 Otherwise, there exists a subsequence of $\{f_n\}$ (that we continue to call $\{f_n\}$) such that $f_n(z)\neq 0$ in $\Delta(0,\delta)$.
 Since $f_n(z)\neq h(z)$, we have $f_n$ and $h$ have no common poles for each $n$ and $f_n(z)-h(z)\neq 0$.
 By Lemma \ref{lschwick}, $\{f_n\}$ is normal at $0$. A contradiction.

 Hence, taking a subsequence and renumbering if necessary,
 we may assume that $b_n$ be the zero of $\{f_n\}$ of smallest modulus. Obviously we have $b_n\rightarrow 0$ as $n\rightarrow\infty$. Since $f_n(0)\neq 0$, we have $b_n\neq 0$.
 Let $G_n(\zeta)=b_n^{k-1}f_n(b_n\zeta)$, we have
\begin{lista}
\item [$(k1)$] $G_n(\zeta)\neq 0$ in $\Delta$; \item [$(k2)$] all
zeros of $G_n(\zeta)$ are multiple; \item [$(k3)$]
$G'_n(\zeta)\neq \frac{\widehat h(b_n\zeta)}{\zeta^k}$ and
$G_n(1)=0$.
\end{lista}
By Lemma \ref{lnormalzz}, $\{G_n(\zeta)\}$ is normal at $\Delta$.
By Lemma \ref{lquasichiefb}, $\{G_n(\zeta)\}$ is quasinormal in
$\mathbb{C}$. Thus, there exists a subsequence of $\{G_n(\zeta)\}$
(that we continue to call $\{G_n(\zeta)\}$) and $E_5 \subset
\mathbb{C}$ such that
\begin{lista}
\item [$(l1)$] $E_5$ have no accumulation point in $\mathbb{C}$;
\item [$(l2)$] $G_n(\zeta)\overset{\chi}{\Longrightarrow}G(\zeta)$
on $\mathbb{C}\backslash E_5$; \item [$(l3)$] for each $\zeta_0\in
E_5$, no subsequence of $\{G_n(\zeta)\}$ is normal at $\zeta_0$.
\end{lista}
Obviously, $E_5\cap\Delta=\emptyset$ and all zeros of $G(\zeta)$
are multiple in $\mathbb{C}\backslash E_5$.

{\bf Case 3.2.1:} \textit{$E_5$ is the empty set.}

Since $G_n(1)=0$, then $G(1)=0$ and hence $1$ is a multiple zero
of $G(\zeta)$ and $G(\zeta)$ is a meromorphic function. Since
$G'_n(\zeta)\neq \frac{\widehat h(b_n\zeta)}{\zeta^k}$ and
$\frac{\widehat
h(b_n\zeta)}{\zeta^k}\overset{\chi}{\Longrightarrow}\frac{1}{\zeta^k}$
in $\mathbb{C}$, then by Lemma \ref{lhurwitz}, we have either
$G'(\zeta)\equiv \frac{1}{\zeta^k}$ in $\mathbb{C}$, or
$G'(\zeta)\neq\frac{1}{\zeta^k}$ in $\mathbb{C}$.

 If $G'(\zeta)\equiv\frac{1}{\zeta^k}$ in $\mathbb{C}$,
 then we have $G(\zeta)=\frac{1}{1-k}(\frac{1}{\zeta^{k-1}}+c)$.
 Obviously, $1$ is not a multiple zero of $G(\zeta)$, a contradiction.

If $G'(\zeta)\neq\frac{1}{\zeta^k}$ in $\mathbb{C}$, then by
Theorem \ref{trational} and Lemma \ref{lrationalb}, we have
$G(\zeta)=c$. Since $G(1)=0$, we have $G(\zeta)\equiv 0$. i.e.,
\begin{equation}\label{equation321a}
G_n(\zeta)=b_n^{k-1}f_n(b_n\zeta)\Rightarrow 0, \quad\text{in}
\quad\mathbb{C}.
\end{equation}
Since $G_n(\frac{a_n}{b_n})=b_n^{k-1}f_n(a_n)=\infty$, by
Hurwitz's Theorem and (\ref{equation321a}), we have
$\frac{a_n}{b_n}\rightarrow \infty$. Writing
$r_n=\frac{b_n}{a_n}$. Hence $r_n\rightarrow 0$ as $n\rightarrow
\infty$.

Since $F_n(r_n)=a_n^{k-1}f_n(b_n)=0$, we have $F(0)=0$ and hence
$F(z)$ is a meromorphic function in $\mathbb C\setminus E_4$.
Since $F'_n(\zeta)\neq \frac{\widehat h(a_n\zeta)}{\zeta^k}$ and
$\frac{\widehat
h(a_n\zeta)}{\zeta^k}\overset{\chi}{\Longrightarrow}\frac{1}{\zeta^k}$
in $\mathbb{C}$, we have by Lemma \ref{lhurwitz}  either
$F'(\zeta)\equiv\frac{1}{\zeta^k}$ in $\mathbb C\setminus E_4$, or
$F'(\zeta)\neq\frac{1}{\zeta^k}$ in $\mathbb C\setminus E_4$.

 Now, we consider the following two cases.

\textbf{Case 1:} \textit{$E_4$ is the empty set.}

If $F'(\zeta)-\frac{1}{\zeta^k}\equiv 0$ in $\mathbb{C}$, then we
have $F(0)=\infty$ which contradicts that $F(0)=0$. If
$F'(\zeta)\neq\frac{1}{\zeta^k}$ in $\mathbb{C}$, then by Theorem
\ref{trational} and Lemma~\ref{lrationalb}, we have $F(\zeta)=c$
which contradicts that $F(1)=\infty$.

\textbf{Case 2:} \textit{$E_4$ is not the empty set.}

 Let $\zeta_0\in E$.
Obviously, $\zeta_0\neq 0$. By Lemma \ref{lquasichiefa}, we have
$F(\zeta)=\int_{\zeta_0}^\zeta\frac{1}{\xi^k}d\xi=\frac{1}{1-k}\big(\frac{1}{\zeta^{k-1}}-\frac{1}{{\zeta_0}^{k-1}}\big)$.
Obviously, we have $F(0)=\infty$ which contradicts $F(0)=0$.

{\bf Case 3.2.2:} \textit{$E_5$ is not the empty set.}

Let $\zeta_0\in E_5$. Obviously, $\zeta_0\neq 0$.
By Lemma \ref{lquasichiefa}, we have
$$G(\zeta)=\int_{\zeta_0}^\zeta\frac{1}{\xi^k}d\xi=\frac{1}{(k-1)\zeta_0^{k-1}}\left(\frac{\zeta^{k-1}-\zeta_0^{k-1}}{\zeta^{k-1}}\right),\quad \zeta \in \mathbb{C}\backslash E_5.$$
Let $\zeta_i$ $(i=1,2,\ldots,k-1)$ be the $k-1$ roots of the equation
 $\zeta^{k-1}-\zeta_0^{k-1}=0$.

Do as in  Case 2.2, we have
$E_5=\{\zeta_1,\zeta_2,\ldots,\zeta_{k-1}\}$.

We claim that $1\in E_5$ and hence
$G(\zeta)=\frac{1}{1-k}\big(\frac{1}{\zeta^{k-1}}-1\big)$. Indeed,
if $1\not\in E_5$, i.e.,  $1-\zeta_0^{k-1}\neq 0$, then
$G(1)=\frac{\zeta_0^{k-1}-1}{(1-k)\zeta_0^{k-1}}\neq 0$. On the
other hand, since $G_n(1)=0$,   we have $G(1)=0$, a contradiction.

Now we have
\begin{equation}\label{equation322a}
   G_n(\zeta)\overset{\chi}{\Longrightarrow}\frac{\zeta^{k-1}-1}{(k-1)\zeta^{k-1}},\quad \zeta\in \Delta.
\end{equation}

By Hurwitz's Theorem, there exist $\gamma_{n,i}$,
$i=1,2,\ldots,k-1$
 such that $\gamma_{n,i}\rightarrow 0$ and $G_n(\gamma_{n,i})=\infty$.
 Since $f$ and $h$ have no common poles, then $G_n(0)\neq \infty$, and we have that  $\gamma_{n,i}\neq 0$ $(i=1,2,\cdots,k-1)$.

Let $s_{n}$ be  one of
$\{\gamma_{n,1},\gamma_{n,2},\cdots,\gamma_{n,k-1}\}$ of largest
modulus. Set $U_n(\xi)=s_n^{k-1}G_n(s_n\xi)$. Obviously, for $n$
sufficiently large, $U_n(\xi)$ have only $k-1$ poles
$\eta_{n,i}=\frac{\gamma_{n,i}}{s_n}$ on $\overline \Delta$. By
(\ref{equation322a}), we have for each $R>0$,
$$U_n(\xi) \neq 0\quad \xi\in\Delta(0,R)$$ for $n$ sufficiently large.
By $(k_3)$ $U'_n(\xi)\neq \frac{\widehat h(b_ns_n\xi)}{\xi^k}$,
and then we get by Lemma \ref{lnormalzz} that $U_n(\xi)$ is normal
in $\mathbb{C}$. Without loss of generality, we assume that
$U_n(\xi)\overset{\chi}{\Longrightarrow}U(\xi)$ in $\mathbb{C}$
and $\eta_{n,i}\rightarrow \eta_i$. Since $U_n(1)=\infty$, we have
$U(1)=\infty$.

We claim that $U(\xi)\equiv \infty$ in $\mathbb{C}$. Otherwise, by
Lemma \ref{lhurwitz},
 either $U'(\zeta)\equiv\frac{1}{\xi^k}$ in $\mathbb{C}$, or $U'(\xi)\neq\frac{1}{\xi^k}$ in $\mathbb{C}$.
If $U'(\xi)\equiv\frac{1}{\xi^k}$ in $\mathbb{C}$, then we have
$U(\xi)=\frac{1}{1-k}(\frac{1}{\xi^{k-1}}+c)$ which contradicts
that $U(1)=\infty$. If $U'(\xi)\neq\frac{1}{\xi^k}$ on
$\mathbb{C}$, then, by Theorem \ref{trational} and Lemma
\ref{lrationalb}, we have $U(\xi)=c$ which contradicts
$U(1)=\infty$.

Now, we have
$$U_n(\xi)\overset{\chi}{\Longrightarrow}\infty,\quad \xi\in \mathbb{C}.$$
Let
$U^\ast_n(\xi)=U_n(\xi)\cdot\prod\limits_{i=1}^{k-1}(\xi-\eta_{n,i})$.
By the maximum principle applied to $/U_n^*$, we get that
$$U^\ast_n(\xi)\overset{\chi}{\Longrightarrow}\infty,\quad \xi\in \mathbb{C}.$$
Let
$G_n^\ast(\zeta)=G_n(\zeta)\cdot\prod\limits_{i=1}^{k-1}(\zeta-\gamma_{n,i})=G_n(\zeta)\cdot\prod\limits_{i=1}^{k-1}(\zeta-s_n\eta_{n,i})$.
By (\ref{equation322a}), for sufficiently large $n$,
$G_n^\ast(\zeta)$ have no pole in $\Delta(0,\frac{1}{2})$ and
hence, by the maximum principle,
$$G_n^\ast(\zeta)\overset{\chi}{\Longrightarrow}\frac{\zeta^{k-1}-1}{k-1},\quad \zeta\in \Delta(0,\frac{1}{2}).$$
So we have $G_n^\ast(0)\rightarrow \frac{1}{1-k}$ as $\rightarrow\infty$. On the other hand,
\begin{align}\nonumber\label{equation322b}
  G_n^\ast(s_n\xi)&=G_n(s_n\xi)\cdot\prod\limits_{i=1}^{k-1}(s_n\xi-\gamma_{n,i})=G_n(s_n\xi)\cdot\prod\limits_{i=1}^{k-1}(s_n\xi-s_n\eta_{n,i})\\
  &=s_n^{k-1}G_n(s_n\xi)\cdot\prod\limits_{i=1}^{k-1}(\xi-\eta_{n,i})=U_n^\ast(\xi)\overset{\chi}{\Longrightarrow}\infty,\quad \xi\in \mathbb{C}.
\end{align}
(\ref{equation322b}) implies that $G_n^\ast(0)\rightarrow \infty$
as $\rightarrow\infty$ which contradicts $G_n^\ast(0)\rightarrow
\frac{1}{1-k}$ as $\rightarrow\infty$. The proof of the Theorem is
complete.
\end{proof}

\end{document}